\documentclass [11pt]{amsart} 
\usepackage{epic,eepic,latexsym, amssymb, amscd, amsfonts, color}

\input xy
\xyoption {all}
\def \mv {{\mathfrak M^{+}_{v}}}  

\def \mw {{\mathfrak M^{+}_{w}}}

\def \gbar {{\bar g}}
\def \rp {{{\bf R}}}
\def \fm {{{\bf R}{\mathcal S}}}
\def \fmd {{{\bf R}{\mathcal S}^{\dagger}}}

\newtheorem{theorem}{Theorem} 
\newtheorem {lemma}{Lemma} 
 
\newtheorem {proposition}{Proposition}

\newcommand{\comment}[1]{}
\newtheorem*{thma}{{Theorem}}
\newtheorem*{lemmaa}{{Lemma 1A}}
\newtheorem*{propositiona}{{Proposition 1A}}
\newtheorem*{th1a}{{Theorem 2A}}
\newtheorem* {th1r}{{Theorem 2R}}

\theoremstyle{definition}
\newtheorem{remark}{Remark}

\begin{document}

\baselineskip=17pt
\title[On the strange duality conjecture for abelian surfaces]{On the strange duality conjecture for abelian surfaces}
\author {Alina Marian}
\address {Department of Mathematics}
\address {Northeastern University}
\email {a.marian@neu.edu}
\author {Dragos Oprea}
\address {Department of Mathematics}
\address {University of California, San Diego}
\email {doprea@math.ucsd.edu}

\begin{abstract}
We study Le Potier's strange duality conjecture for moduli spaces of sheaves over generic abelian surfaces. We prove the isomorphism for abelian surfaces which are products of elliptic curves, when the moduli spaces consist of sheaves of equal ranks and fiber degree $1$. The birational type of the moduli space of sheaves is also investigated. Generalizations to arbitrary product elliptic surfaces are given. \end {abstract}
\maketitle
\section{Introduction}

There are three versions of Le Potier's strange duality conjecture for abelian surfaces, formulated and supported numerically in \cite {abelian}. In this article, we confirm two of them for product abelian surfaces. This establishes the conjecture over an open subset in the moduli space of polarized abelian surfaces. 

\subsection{The strange duality morphism} To set the stage, let $(X,H)$ be a polarized complex abelian surface. For a coherent sheaf $V\to X$, denote by $$v = \text{ch} \, V \in H^{\star} (X, {\mathbb Z})$$ its Mukai vector. Fix two Mukai vectors $v$ and $w$ such that the orthogonality condition $$\chi(v\cdot w)=0$$ holds. We consider the two moduli spaces $\mathfrak M_v^{+}$ and $\mathfrak M_w^{+}$ of $H$-semistable sheaves of type $v$ and $w$ with fixed determinant. They carry two determinant line bundles $$\Theta_w\to \mathfrak M_v^{+}, \,\,\, \Theta_v\to \mathfrak M_w^{+},$$ whose global sections we seek to relate.  

More precisely, according to Le Potier's strange duality conjecture \cite{lepotier} \cite {abelian}, if the condition $$c_1(v \cdot w) \cdot H > 0$$ is satisfied, the jumping locus $$\Theta^+_{vw}=\{(V, W): h^1(V\otimes W)\neq 0\}\subset \mathfrak M_v^{+}\times \mathfrak M_w^{+}$$ is expected to be a divisor. Further, its defining equation is a section of the split bundle $$\Theta_w\boxtimes \Theta_v\to \mv\times \mw,$$ conjecturally inducing an isomorphism $$\mathsf D^{+}:H^0(\mathfrak M_v^{+}, \Theta_w)^{\vee}\to H^0(\mathfrak M_w^{+}, \Theta_v).$$ 

\noindent {This expectation is supported numerically: it was established in \cite {abelian} that $$\chi (\mathfrak M_v^+, \Theta_w)=\chi (\mathfrak M_w^+, \Theta_v)=\frac{\chi (X, L^+)}{d_v+d_w} \binom{d_v+d_w}{d_v},$$ where $L^+$ is a line bundle on $X$ with $$c_1 (L^+)=c_1 (v\cdot w)$$ and $$d_v=\frac{1}{2} \dim \mathfrak M_v^+,\,\,\, d_w=\frac{1}{2} \dim \mathfrak M_w^+.$$}

Similarly, letting $\mathfrak M_v^-$ and $\mathfrak M_w^-$ denote the moduli spaces of sheaves with fixed determinant of their Fourier-Mukai transforms, we have the symmetry $$\chi  (\mathfrak M_v^-, \Theta_w)=\chi (\mathfrak M_w^-, \Theta_v)=\frac{\chi (X, L^-)}{d_v+d_w} \binom{d_v+d_w}{d_v},$$ where $L^-$ is a line bundle on $X$ with $$c_1 (L^-)=c_1 (\widehat v\cdot \widehat w),$$ the hats denoting the Fourier-Mukai transforms. The map  $$\mathsf D^-:H^0(\mathfrak M_v^{-}, \Theta_w)^{\vee}\to H^0(\mathfrak M_w^{-}, \Theta_v)$$ induced by the theta divisor $\Theta^-_{vw}\subset \mathfrak M_v^{-}\times \mathfrak M_w^{-}$ is expected to be an isomorphism. 

\subsection {Results} We will establish the isomorphisms for abelian surfaces which split as products of elliptic curves $$X=B\times F.$$ We regard $X$ as a trivial fibration $\pi_B: X\to B$, and write $\sigma$ and $f$ for the class of the zero section and of the fiber over zero. We assume that the polarization $H$ is suitable in the sense of \cite {F} {\it{i.e.}} $$H=\sigma+ Nf, \text{ for } N\gg 0.$$ 

Over {\it simply connected} elliptic surfaces, for coprime rank and fiber degree, the moduli space of sheaves is birational to the Hilbert scheme of points, as shown by Bridgeland \cite {B}. For sheaves with fixed determinant, the situation is subtler over elliptic {\it abelian} surfaces. A refinement of Bridgeland's argument, using a Fourier-Mukai transform with kernel given by a universal Atiyah bundle, allows us to prove the following result.
\begin {theorem} \label{t3}Let $X= B\times F$ be a product abelian surface, and let $v$ be a Mukai vector such that the rank $r$ and the fiber degree $d=c_1(v)\cdot f$ are coprime. Then, the moduli space $\mathfrak M_v^+$ is birational to 
$$\mathfrak X_v^{+}=\{(Z, b):  a_B(Z)=r b\}\subset X^{[d_v]}\times B,$$ where $a_B$ denotes the addition over the base elliptic curve $B$. 
\end {theorem}
\noindent A similar statement holds for the moduli space $\mathfrak M_v^{-}.$ By contrast, the generalized Kummer variety $K_v$ associated to the higher rank vectors is birational to the Kummer variety in rank $1$, as noted in \cite{Y}. This fact is recovered in two ways
while establishing Theorem \ref{t3}.

The following two theorems capture our main results concerning strange duality. 
\begin {theorem} \label{t1}Let $X= B\times F$ be a product abelian surface. Let $v$ and $w$ be two orthogonal Mukai vectors of equal ranks $r\geq 3$, with $$c_1(v)\cdot f=c_1(w)\cdot f=1.$$
Then, $$\mathsf D^+:H^0(\mathfrak M^+_v, \Theta_w)^{\vee}\to H^0(\mathfrak M^+_w, \Theta_v)$$ is an isomorphism. 

\end{theorem}

Similarly, we show

\begin {theorem}\label{t2} Let $X=B\times F$ be a product abelian surface. Assume $v$ and $w$ are two orthogonal Mukai vectors of ranks $r, s\geq 3$ and equal Euler characteristics $\chi=\chi'$, with $$c_1(v)\cdot f=c_1(w)\cdot f=1.$$ 
Then, $$\mathsf D^-:H^0(\mathfrak M^-_v, \Theta_w)^{\vee}\to H^0(\mathfrak M^-_w, \Theta_v)$$ is an isomorphism. 
\end {theorem}

In particular, the theorems imply that $\Theta^{\pm}$ are divisors in the products $\mathfrak M_v^{\pm}\times \mathfrak M_w^{\pm}.$ 
\vskip.1in

\subsubsection{Higher genus} The requirement that $B$ be elliptic can in fact be removed in Theorem \ref{t1}. Indeed, consider $C$ a smooth projective curve of genus $g\geq 1$, and let $\sigma$ denote the zero section of the trivial fibration $$X=C\times F\to C.$$ We show

\begin {th1a} \label{t4} Assume that $X=C\times F$ is a product surface as above, and let $v, w$ be orthogonal Mukai vectors of equal ranks $r\geq g$ and $r\neq 2$, such that 
\begin {itemize}
\item [(i)] the determinants are fixed of the form 
$$\det v=\mathcal O(\sigma)\otimes \ell_v,\,\,\,\det w=\mathcal O(\sigma)\otimes \ell_w,$$ for generic line bundles $\ell_v$ and $\ell_w$ of fixed degree over the curve $C$;
\item [(ii)] $\dim \,\mv+\dim\, \mw\geq 8r(g-1).$\end{itemize}\noindent Then $$\mathsf D^+:H^0(\mathfrak M^+_v, \Theta_w)^{\vee}\to H^0(\mathfrak M^+_w, \Theta_v)$$ is an isomorphism. 

\end {th1a}

\subsection {Comparison} Theorems \ref{t1} and \ref{t2} parallel the strange duality results for {\it simply connected} elliptic surfaces. Let $$\pi:Y\to \mathbb P^1$$ be a simply connected elliptic fibration with a section and at worst irreducible nodal fibers. The dimension of the two complementary moduli spaces $\mathfrak M_v$ and $\mathfrak M_w$ will be taken large enough compared to the constant $$\Delta= \chi(Y, \mathcal O_Y) \cdot \left((r+s)^2+(r+s)+2\right)-2(r+s).$$ The polarization is still assumed suitable. The following was proved by combining \cite {generic} and \cite {BH}:
\begin {thma} \label{ellsur} Let $v$ and $w$ be two orthogonal topological types of rank $r, s\geq 3$, such that \begin {itemize} \item [(i)] the fiber degrees $c_1(v)\cdot f=c_1(w)\cdot f=1,$
\item [(ii)] $\text {dim }\mathfrak M_v+\text{dim } \mathfrak M_w\geq \Delta.$ \end {itemize} Then, $$\mathsf D:H^0(\mathfrak M_v, \Theta_w)^{\vee}\to H^0(\mathfrak M_w, \Theta_v)$$ is an isomorphism. 
\end {thma}

We suspect that an analogous statement can be made for all (not necessarily simply connected, with possibly reducible fibers) elliptic surfaces $Y\to C$, going beyond the scope of Theorems \ref{t1} and 2A. 

The results for simply connected fibrations and abelian surfaces both rely on Fourier-Mukai techniques, but the geometry is more involved in the abelian case, as already illustrated by the birationality statement of Theorem \ref{t3}. Via Fourier-Mukai, instead of a rather standard analysis of tautological line bundles over Hilbert scheme of points in the simply connected case, one is led here to studying sections of suitable theta bundles over the schemes $\mathfrak X_v^+$. This requires new ideas. We prove the duality for spaces of sheaves of equal ranks and fiber degree $1$, but believe these assumptions may be relaxed. The case $r=s=2$ is also left out of our theorems: while the Fourier-Mukai arguments do not cover it, we believe strange duality holds here as well.

\subsection {Variation in moduli} Let $\mathcal A_d$ be the moduli space of pairs $(X, H)$, where $X$ is an abelian surface and $H$ is an ample line bundle inducing a polarization of type $(1,d)$ on $X$. For a Mukai vector $v$ with $$c_1(v) = c_1 (H),$$ we consider the relative moduli space of sheaves $$\pi: {\mathfrak M} [v]^{+} \longrightarrow \mathcal A_d$$ whose fiber over a surface $(X, H)$ is the moduli space of $H$-semistable sheaves with Mukai vector $v$ and fixed determinant equal to $H$.  
Associated with two orthogonal Mukai vectors $v$ and $w$, there is a universal canonical Theta divisor $$\Theta^+_{vw} = \{(X, H, V, W): \, h^1 (X, V \otimes W) \neq 0 \} \subset {\mathfrak M}[v]^{+} \times_{\mathcal A_d} {\mathfrak M} [w]^{+}, $$ giving rise to a line bundle which splits as a product $$\Theta^+_{vw} = \Theta_w \boxtimes \Theta_v \, \, \, \text{on} \, \, \, {\mathfrak M}[v]^{+} \times_{{\mathcal A}_{d}} {\mathfrak M} [w]^{+}. $$  Pushforward via the morphisms $${\mathfrak M}[v]^{+}\to \mathcal A_d,\,\,\,\, {\mathfrak M}[w]^{+}\to \mathcal A_d$$ yields two coherent sheaves of generalized theta functions over $\mathcal A_d$, $${\mathbb W} = R^0\pi_{\star} \Theta_w , \, \, \, \,   {\mathbb V} = R^0 \pi_{\star} \Theta_v .$$ Let ${\mathcal H} \subset \mathcal A_d$ be the Humbert surface parametrizing split abelian surfaces $(X, H)$ with $$X = B \times F\text{ and }H = L_B \boxtimes L_F,$$ for line bundles $L_B, L_F$ over $B$ and $F$ of degrees $d$ and $1$ respectively. 

Theorem \ref{t1} can be rephrased as the following generic strange duality statement:
\begin{th1r} Assume $v$ and $w$ are orthogonal Mukai vectors of equal ranks $r\geq 3$ with \begin {itemize}
\item [(i)] $c_1(v)=c_1(w)=H;$
\item [(ii)] $\langle v, v\rangle \geq 2(r^2+r-1), \,\, \langle w, w\rangle \geq 2(r^2+r-1).$
\footnote {Assumption (ii) allows us to exchange $H$-stability with stability with respect to a suitable polarization, since in this case the ensuing moduli spaces agree in codimension $1$. This was proved for $K3$ surfaces in the appendix of \cite {generic}. The case of abelian surfaces follows by the same argument.}
\end {itemize}
Then, the sheaves ${\mathbb V}$ and ${\mathbb W}$ are locally free when restricted to the Humbert surface ${\mathcal H}$, and $\Theta^+_{vw}$ induces an isomorphism $${\mathsf D}: {\mathbb W}^{\vee} \longrightarrow {\mathbb V} \, \, \, \, \text{along} \, \, \, {\mathcal H}. $$ 
\end{th1r}

\subsection {Acknowledgements} The authors acknowledge support from the NSF via grants DMS 1001604 and DMS 1001486, and thank CRM Barcelona, MPI Bonn, and the University of Bonn for their hospitality in the summer of 2012. 

\section{Moduli spaces of sheaves of fiber degree $1$} 
\label{s2}

\subsection{Setting}  \label{basiccounts} We consider a complex abelian surface $X$ which is a product of two elliptic curves,
$$X \simeq B \times F.$$ 
Letting $o_B, o_F$ denote the origins on $B$ and $F$, we write $\sigma$ and $f$ for the divisors $B \times o_F$ respectively $o_B \times F$ in the product $B\times F$. Note that $$\sigma^2=f^2=0,\,\,\, \sigma\cdot f=1.$$ 
\subsubsection {Line bundles over $X$} For any positive integers $a$ and $b$, the line bundle $\mathcal O (a\sigma + b f)$ on $X$ is ample. Its higher cohomology vanishes, and we have 
\begin{equation}
\label{sectioncount1}
h^0 (X, {\mathcal O} (a \sigma + b f)) = \chi  (X, {\mathcal O} (a \sigma + b f)) = ab.
\end{equation} Letting $$\pi_B:X\to B, \, \pi_F:X\to F$$ be the natural projections from $X$, we also note that for any $\ell > 0,$
\begin{eqnarray}
\label{sectioncount2}
h^0 (X, {\mathcal O} (\ell\sigma)) &= &h^0 (X, \pi_F^{\star} {\mathcal O} (\ell o_F)) = h^0 (F, {\mathcal O} (\ell o_F)) = \ell, \\  \nonumber
h^0 (X, {\mathcal O} (\ell f)) &=& h^0 (X, \pi_B^{\star} {\mathcal O} (\ell o_B)) = h^0 (B, {\mathcal O} (\ell o_B)) = \ell.
\end{eqnarray}
Thus every section of ${\mathcal O}_X (\ell f)$ corresponds to a divisor of the form $$\bigsqcup_{i=1}^{\ell} z_i \times F, \, \, \, \text{for} \, \, \, z_1 + \cdots + z_{\ell} = o_B.$$
Furthermore, equations  \eqref{sectioncount1} and \eqref{sectioncount2} imply that every section of the line bundle ${\mathcal O}_X (\sigma + \ell f)$ vanishes along a divisor of the form
$$\sigma \cup \bigsqcup_{i=1}^{\ell} z_i \times F, \, \, \, \text{for} \, \, \, z_1 + \cdots + z_{\ell} = o_B.$$

\subsubsection{Fourier-Mukai functors}
Two different Fourier-Mukai transforms will be considered on $X$. First, there is the Fourier-Mukai transform with respect to the standardly normalized Poincar\'{e} bundle $$\mathcal P\to X\times \widehat X.$$ We identify $\widehat X\cong X$ using the polarization $o_B\times F+B\times o_F$. Specifically, our (perhaps non-standard) sign choice is so that $y=(y_B, y_F) \in X$ is viewed as a degree 0 line bundle $y \to X$ via the association $$y \mapsto {\mathcal O} _{B} (y_B - o_B)\boxtimes \mathcal O_F(y_F-o_F).$$ The Poincar\'{e} bundle is given then by $$ \mathcal P \to X \times X, \, \, \, \ \mathcal {P} =\mathcal P_B\boxtimes \mathcal P_F$$ where $$\mathcal P_B\to B\times B,\,\,\, \mathcal P_B= {\mathcal O}_{B \times B} (\Delta_B)\otimes p^{\star} {\mathcal O}_B (-o_B) \otimes q^{\star} {\mathcal O}_B (-o_B), \text{ and } $$ $$\mathcal P_F\to F\times F, \, \, \, \, {\mathcal P}_F = {\mathcal O}_{F \times F} (\Delta_F) \otimes p^{\star}{\mathcal O}_F (-o_F) \otimes q^{\star} \mathcal O_F (-o_F).$$
The Fourier-Mukai transform with kernel $\mathcal P$ is denoted $$\fm:\mathbf D(X)\to \mathbf D( X).$$ 

A second Fourier-Mukai transform is defined by considering the relative Picard variety of $\pi_B: X \to B.$ We identify $F \cong \widehat F$ so that $$ y_F \in F \mapsto {\mathcal O}_F (y_F - o_F),$$ and we  let  $$\fm^{\dagger}: \mathbf D(X)\to \mathbf D(X)$$ be the Fourier-Mukai transform whose kernel is the pullback of the Poincar\'e sheaf $$\mathcal P_F\to F\times F$$ to the product $X\times_{B} X \cong F \times F\times B.$  Both Fourier-Mukai transforms are known to be equivalences of derived categories.

\vskip.1in

Several properties of the Fourier-Mukai will be used below. First, for integers $a, b$ we have $$\det \fm(\mathcal O(a\sigma+bf))=\mathcal O(-b\sigma-af),$$ and $$\det \fm^{\dagger} (\mathcal O(a\sigma+bf))=\mathcal O(-\sigma+ab f).$$ Next, we have the following standard result, which will in fact be proved in greater generality in Lemma 1A of Section \ref{arbitraryfiber}. 

\begin {lemma} \label{l1} For $x=(x_B, x_F)\in X,\,\, y=(y_B, y_F)\in \widehat X,$ we have $$\fm^{\dagger}(t_{x}^{\star} E)=t_{x_B}^{\star} \fm^{\dagger}(E)\otimes x_F^{\vee},$$ $$\fm^{\dagger}(E\otimes y)=t_{y_F}^{\star}\fm^{\dagger} (E)\otimes y_B.$$
\end {lemma}

\subsection{Moduli spaces of sheaves}

We  consider sheaves over $X$ of Mukai vector $v$ such that $$\text{ rank }v=r,\,\,\, \chi(v)=\chi,\,\,\, c_1(v)\cdot f=1.$$ We set $$d_v=\frac{1}{2}\langle v, v\rangle = \frac{c_1(v)^2}{2}-r\chi,$$ which is half the dimension of the moduli space $\mathfrak M_v^+$ below. Recall from the introduction that the polarization $H$ is suitable {\it i.e.}, 
$$H = \sigma + N f, \, \,\, \text{for}\, \,\, N \gg 0.$$  We are concerned with three moduli spaces of sheaves: 
\begin{itemize}
\item[(i)] 
The moduli space $K_v$ of $H$-semistable sheaves $V$ on $X$ with $$\det V=\mathcal O(\sigma+m f),\,\,\, \det \fm(V)=\mathcal O(-m \sigma-f).$$ Thus the determinants of the sheaves and of the Fourier-Mukai transforms of sheaves in $K_v$ are fixed. Here, we wrote $$c_1(v)=\sigma+mf$$ where $$m=d_v+r\chi.$$ 
\item[(ii)]
The moduli space $\mv$ of $H$-semistable sheaves $V$ on $X$ with $$\det V=\mathcal O(\sigma+m f).$$ 
The two spaces (i) and (ii) are related via the degree $d_v^4$ \'{e}tale morphism $$\Phi_v^{+}:K_{v}\times X\to \mathfrak M_v^{+},\, \, \, \, \Phi_v^{+} (V, x) = t_{rx}^{\star} V \otimes t_x^{\star} \det V^{-1} \otimes \det V.$$ In explicit form, we write equivalently
\begin{equation}
\label{psiplus}
\Phi^{+}_v(V, x)=t_{rx}^{\star} V\otimes (x_F\boxtimes x_B^{m})
\end{equation}
\item [(iii)] Finally, there is a moduli space $\mathfrak M_v^{-}$ of semistable sheaves whose Fourier-Mukai transform has fixed determinant $$\det \fm (V)=\mathcal O(-f-m\sigma).$$ In this case, we shall make use of the \'{e}tale morphism $$\Phi_v^{-}:K_v\otimes \widehat X \to \mathfrak M_v^{-}, \, \, \, \, \Phi_v^{-}(V, y)=t_{(y_B, my_F)}^{\star}V\otimes y^{\chi}.$$ Note that we have indeed  \begin{eqnarray*}\det \fm \Phi_v^-(V, y)&=&\det \fm (t_{(y_B, my_F)}^{\star}V\otimes y^{\chi})=t_{\chi y}^{\star} \det \fm (t_{(y_B, my_F)}^{\star}V)\\ &=&t_{\chi y}^{\star} \det (\fm (V)\otimes (y_B, my_F)^{-1})=t_{\chi y}^{\star} \mathcal O(-f-m\sigma)\otimes (y_B, my_F)^{-\chi}\\ &=&\mathcal O(-f-m\sigma).\end{eqnarray*}

\end{itemize} 
For further details regarding the morphisms $\Phi_v^{+}, \, \Phi_v^{-}$, we refer the reader to Sections 4 and 5 of \cite{abelian}.
 
\vskip.2in

\subsection{Birationality of $K_v$ with the generalized Kummer variety $K^{[d_v]}$}\label{bir}

We now establish in two different ways an explicit birational map $$\Psi_r:K^{[d_v]} \dasharrow K_v,$$ where $K^{[d_v]}$ is the generalized Kummer variety of dimension $2d_v -2,$ $$K^{[d_v]}=\{Z\in X^{[d_v]}: \, a(Z)=0\}.$$ As usual, $a:X^{[d_v]}\to X$ denotes the addition map on the Hilbert scheme. 

\subsubsection{O'Grady's construction}\label{sss} We first obtain the map $\Psi_r$ by induction on the rank $r$ of the sheaves, following O'Grady's work \cite{ogrady} for elliptically fibered $K3$ surfaces. 
To start the induction, we let $Z \subset X$ be a zero dimensional subscheme of length $$\ell = \ell(Z) = d_v,$$ such that $a(Z)=0$, and further satisfying the following conditions:
\begin{itemize}
\item [(i)] $Z$ does not contain two points in the same fiber of $\pi_B,$
\item [(ii)] $Z$ does not intersect the section $\sigma.$ 
\end{itemize}
Such a $Z$ is a generic point in the generalized Kummer variety $K^{[d_v]}.$ 
We let $$m_1 = \chi + d_v,$$ which the correct number of fibers needed when $r=1$, and set
$$ V_1=I_Z\otimes \mathcal O(\sigma+m_1 f).$$ 
Note that $\chi (V_1) = \chi,$ and that $V_1$ thus constructed belongs to $K_v$. The requirement $$\det \fm(V_1)=\mathcal O(-m_1\sigma-f),$$ is a consequence of the fact that $a(Z)=0.$

We claim that for $Z$ as above we have 
\begin{equation}
\label{basecase}
h^1(V_1(-\chi f))=h^1(I_Z\otimes \mathcal O(\sigma +\ell f))=h^0(I_Z(\sigma+\ell f))=1.
\end{equation}
Indeed, as explained in Section \ref{basiccounts}, every section of  $\mathcal O_X(\sigma+\ell f)$ vanishes along a divisor of the form $$\sigma+\pi_B^{\star}\, ({z_1}+\ldots {z_\ell)}$$ where $z_1, \ldots, z_{\ell}$ are points of $B$ such that $z_1+\ldots+z_{\ell}=o_B$. 
Conditions (i) and (ii) ensure that there is a unique such divisor passing through $Z$: $z_1, \ldots, z_{\ell}$ are the $B$-coordinates of the distinct points of $Z.$

We inductively construct extensions \begin{equation}\label{ext}0\to \mathcal O(\chi f)\to V_{r+1}\to V_r\to 0,\end{equation} with stable middle term. The sheaves obtained will satisfy $$\chi(V_r(-\chi f))=0.$$ In order to get \eqref{ext}, we show $$\text{ext}^1(V_r, \mathcal O(\chi f))=h^{1}(V_r(-\chi f))=h^0(V_r(-\chi f))=1.$$ For $r=1$, the ext-dimension is 1 by \eqref{basecase}. Assuming the statement for $r$, we get the unique nontrivial extension, and argue for stability of the middle term $V_{r+1}$. Indeed, Proposition I.4.7 of \cite{ogrady} gives stability of the middle term and also asserts the restrictions of \eqref{ext} to any fibers of $\pi$ do not split, if a suitable vanishing hypothesis holds. Precisely, we require that for all fibers $\mathsf f$, we have \begin{equation}\label{vns}\text{Hom}(V_r|_{\mathsf f}, \mathcal O_{\mathsf f})=0.\end{equation} The fiber restrictions are calculated as in Lemma $1$ in \cite {generic} as follows: 
\begin {itemize}
\item [(i)] the restriction of the defining exact sequence to any fiber $f_\eta$ yields $$0\to \mathcal O_{f_{\eta}}\to V_{r+1}|_{f_\eta}\to V_r|_{f_{\eta}}\to 0.$$ Hence inductively, for fibers avoiding $Z$, we have  that \begin {equation}\label{f1}V_{r}|_{f_\eta}\cong {\mathsf E}_{r, o},\end{equation} the Atiyah bundle of rank $r$ and determinant $\mathcal O_{f_\eta}(o)$.
\item [(ii)] for fibers $f_z$ containing points $z\in Z$ we have \begin{equation}\label{f2}V_{r}|_{f_z}=\mathsf E_{r-1, z}\oplus \mathcal O_{f_z}(o-z).\end{equation} \end {itemize} As a consequence of equations \eqref{f1} and \eqref{f2},  the vanishing \eqref{vns} stated above is satisfied.

To complete the argument, we establish now that there exists a unique extension \eqref{ext} for $r+1$ {\it i.e.}, we show $h^1(V_{r+1}(-\chi f))=1.$ First, stability of $V_r$ and $V_{r+1}$ ensures that $$H^2(V_r(-\chi f))=H^2(V_{r+1}(-\chi f))=0,$$ and the cohomology exact sequence of the unique extension for $r$ is $$0\to H^0(\mathcal O)\to H^0(V_{r+1}(-\chi f))\to H^0(V_r(-\chi f))\to H^1(\mathcal O)$$ $$\to H^1(V_{r+1}(-\chi f))\to H^1(V_r(-\chi f))\to H^2(\mathcal O)\to 0.$$ By dimension counting, the last map has to be a bijection, hence $H^1(\mathcal O)\to H^{1}(V_{r+1}(-\chi f))$ is surjective. We claim this map cannot be a bijection. Indeed, if  $H^1(\mathcal O)\to H^{1}(V_{r+1}(-\chi f))$ were bijective, we would have that $$H^0(V_r(-\chi f))\to H^1(\mathcal O)$$ is the zero map. But this map is simply multiplication by the extension class which is nontrivial. Since $h^0(V_{r+1}(-\chi f))\geq 1$ from the same sequence, and $H^0(V_{r+1}(-\chi f))=H^1(V_{r+1}(-\chi f)),$ we conclude both vector spaces are of dimension $1$, as needed. 

Finally, the argument clearly shows that the map $\Psi_r$ is injective. Therefore, we obtain a birational map $$\Psi_r: K^{[d_v]} \dasharrow K_v,$$ by the equality of dimensions for the two irreducible spaces. 

\subsubsection{Fourier-Mukai construction} We next point out that $\Psi_r$ can also be viewed as a fiberwise Fourier-Mukai transform. This is similar to the case of $K3$ surfaces analyzed in \cite{generic}. We show
\begin {proposition}\label{p1} For subschemes $Z$ with $a(Z)=0$, satisfying (i) and (ii), we have $$\fm^{\dagger}(V_r^{\vee})=I_{\widetilde Z}(r\sigma-\chi f)[-1]$$ and $$\fm^{\dagger}(V_r)=I_{Z}^{\vee}(-r\sigma+\chi f),$$ where $\widetilde Z$ is the reflection of $Z$ along the zero section $\sigma$.

\proof We prove the first equality. By the calculation of the restrictions of $V_r$ to the fibers in equations \eqref{f1} and \eqref{f2}, it follows that $V_r^{\vee}$ contains no subbundles of positive degree over each fiber. Therefore, $\fm^{\dagger}(V_r^{\vee})[1]$ is torsion free, by Proposition $3.7$ of \cite{BH}. The agreement of $\fm^{\dagger}(V_r^{\vee})[1]$ and $I_{\widetilde Z}(r\sigma - \chi f)$ holds fiberwise. This can be seen using \eqref{f1} and \eqref{f2},  just as in Proposition $1$ in \cite {generic}. The proof is completed by proving equality of determinants. In turn, this follows inductively from the exact sequence \eqref{ext}: $$\det \fm^{\dagger}(V_{r+1}^{\vee})=\det \fm^{\dagger}(V_r^{\vee})\otimes \det \fm^{\dagger}(\mathcal O(-\chi f))=\det \fm^{\dagger}(V_r^{\vee})\otimes \mathcal O(-\sigma).$$ The base case $r=1$ is immediate, as $$\det \fm^{\dagger} (I_Z^{\vee}(-\sigma-m_1 f))=\det \fm^{\dagger} (\mathcal O(-\sigma - m_1 f))\otimes_{z\in Z}\det (\fm^{\dagger}(\mathcal O_z^{\vee}))^{\vee}$$ $$=\mathcal O(-\sigma + m_1 f)\otimes_{z\in Z} \det (\fm^{\dagger}(\mathcal O_z[-2]))^{\vee}=\mathcal O(-\sigma + m_1 f)\otimes_{z\in Z} \mathcal O(-f_z)$$ $$=\mathcal O(-\sigma+m_1 f)\otimes \mathcal O(-\ell f)=\mathcal O(-\sigma+\chi f),$$ using that $a(Z)=0$. Here $f_z$ denotes the fiber through $z\in Z$. A similar calculation can be carried out for the Fourier-Mukai with kernel $\mathcal P^{\vee}$. 

The second statement follows then by Grothendieck duality.  We refer the reader to Proposition $2$ in \cite {generic} for an identical computation. 

\qed

\end{proposition}
Proposition \ref{p1} gives an explicit description of O'Grady's construction in terms of the Fourier-Mukai transform, at least away from the divisors (i) and (ii). Explicitly, under the correspondence $$K^{[d_v]}\ni Z\mapsto I_{\tilde Z}(r\sigma - \chi f)[-1],$$ the birational map $\Psi_r: K^{[d_v]} \dasharrow K_v$ is given by the inverse of $\fm^{\dagger}$, followed by dualizing $$\Psi_r\cong \mathbf D_X\circ \left(\fm^{\dagger}\right)^{-1}.$$ However, viewed as a Fourier-Mukai transform, $\Psi_r$ extends to an isomorphism in codimension $1$, whenever $r \geq 3.$ This is established in Sections 3 and 5 of \cite{BH} on general grounds, but can also be argued directly. 

Indeed, the inverse of $\fm^{\dagger}$ is, up to shifts, the Fourier-Mukai transform with kernel the dual $\pi_{F\times F}^{\star}{\mathcal P_F}^{\vee}$ of the fiberwise Poincar\'{e} line bundle on $X \times_B X \simeq F \times F\times B.$ We claim that as long as $r \geq 3$ and $Z$ contains no more than two points in the same fiber of $$\pi_B: X \to B,$$ the Fourier-Mukai image of the sheaf $I_Z (r\sigma - \chi f)$ is a vector bundle. Stability is automatic, as the restrictions to all elliptic fibers which do not pass through $Z$ are isomorphic to the Atiyah bundle of rank $r$ and degree $1$, therefore the restriction to the generic fiber is stable. The locus of $Z \in K^{[d_v]}$ with at least $3$ points in the same elliptic fiber has codimension $2.$ 

To see the claim above, note that the restriction of the ideal sheaf of a point to the elliptic fiber $\mathsf f$ through that point is $$I_p{|_{\mathsf f}} = {\mathcal O}_{\mathsf f} (-p) + {\mathcal O}_p,$$ and similarly  $$I_{p, p'}{|_{\mathsf f}} = {\mathcal O}_{\mathsf f} (-p-p') + {\mathcal O}_p+\mathcal O_{p'},$$ for (possibly coincident) points $p, p'$ in the fiber $\mathsf f$. 
Thus, the restriction of $I_Z(r\sigma-\chi f)$ to fibers $\mathsf f$ can take the form $${\mathcal O}_{\mathsf f} (r o_F), \,\,\, {\mathcal O}_{\mathsf f} (r o_F) \otimes \left ( {\mathcal O}_{\mathsf f}(-p) + {\mathcal O}_p \right ) \text{ or }{\mathcal O}_{\mathsf f} (r o_F) \otimes \left ( {\mathcal O}_{\mathsf f}(-p - p') + {\mathcal O}_p + {\mathcal O}_{p'} \right )$$ which are all $IT_0$ with respect to $\mathcal P_F^{\vee},$ for $r\geq 3$. Thus the Fourier-Mukai transform $(\fm^{\dagger})^{-1}$ of $I_Z(r\sigma-\chi f)[-1]$ is a vector bundle by cohomology and base-change. 

When $r=2$, $\Psi_2$ is defined away from the divisor of subschemes $Z \in K^{[d_v]}$ with at least two points in the same elliptic fiber. As $K^{[d_v]}$ and  $K_v$ are irreducible holomorphic symplectic, $\Psi_2$ extends anyway to a birational map which is regular outside of codimension 2, but this extension is no longer identical to the Fourier-Mukai transform. In fact, semistable reduction is necessary to construct the extension, as in Section $I.4$ of \cite {ogrady}. We will not pursue it in this paper.

\subsection{The moduli space $\mv$ via Fourier-Mukai}

We investigate how sheaves of fixed determinant change under Fourier-Mukai. Recall the morphism \eqref{psiplus}, and let $$V= \Phi_v^{+} (E,\, x),$$ for a pair $(E, \, x) \in K_{v}\times X.$ Using Lemma \ref{l1} and Proposition \ref{p1}, we calculate \begin{eqnarray*}\fm^{\dagger}(V)&=&\fm^{\dagger}(t_{rx}^{\star} E\otimes \left(x_B^m\boxtimes x_F\right))=t_{x_F}^{\star} \fm^{\dagger}(t_{rx}^{\star}E)\otimes x_B^{m}\\&=&t_{x_F}^{\star} \left(t_{rx_B}^{\star}\fm^{\dagger}(E)\otimes x_F^{-r}\right)\otimes x_B^{m}\\&=&t_{x_F+rx_B}^{\star}\left( I_{Z}^{\vee} (-r\sigma+\chi f)\right)\otimes \left(x_F^{-r}\boxtimes x_B^{m}\right)\\&=&I_{{Z^{+}}}^{\vee}(-r\sigma+\chi f)\otimes x_B^{d_v},\end{eqnarray*} where $$Z^{+}=t_{x_F+rx_B}^{\star} Z, \, \, \, \text{so that } a_B(Z^{+})=-d_v rx_B.$$  In a similar fashion, we prove that $$\fm^\dagger (V^{\vee})=I_{Z_{+}}(r\sigma - \chi f)[-1]\otimes x_B^{-d_v},$$ where now $$Z_+=t^{\star}_{-x_F+rx_B}\widetilde Z=\widetilde {Z^+}.$$

From the first equation, we obtain the rational map $$\fm^{\dagger}:K_v\times X\dasharrow \mathfrak X^{+}_v$$ where $$\mathfrak X^{+}_v=\left\{(Z^{+}, z_B):  a_B(Z^{+})= r z_B\right\}\subset X^{[d_v]}\times B$$ via the assignment $$(E, x)\mapsto (Z^+, -d_v x_B).$$ This map has degree $d_v^4$ and descends to $\mathfrak M_v^+$. Since $\Phi_v^+$ is also of degree $d_v^4$ and $\mathfrak X_v^+$ is irreducible of the same dimension as $\mathfrak M_v^+$, we obtain a birational isomorphism $$\mathfrak M_v^+\dasharrow \mathfrak X_v^+$$ in such a fashion that 
$$\fm^{\dagger}(V)=I_{Z^+}^{\vee}(-r\sigma+\chi f) \otimes z_B^{-1}$$ and $$\fm^{\dagger}(V^{\vee})=I_{\widetilde {Z^+}}(r\sigma - \chi f)\otimes z_B[-1].$$ The discussion of the previous subsection shows the birational isomorphism is given (explicitly as a Fourier-Mukai transform) away from codimension $2$. 

\vskip.1in

\subsection {The moduli space $\mathfrak M_v^{-}$ via Fourier-Mukai} 

A similar argument applies to the moduli space $\mathfrak M_v^{-}$ of sheaves with fixed determinant of their Fourier-Mukai transform $$\det \fm (V)=\mathcal O(-f-m\sigma).$$ In this case, we have a morphism $$\Phi_v^{-}:K_v\otimes \widehat X \to \mathfrak M_v^{-}$$ given by $$\Phi_v^{-}(E, y)=t_{(y_B, my_F)}^{\star}E\otimes y^{\chi}.$$ We calculate \begin{eqnarray*}\fmd (\Phi_v^{-}(E, y))&=&\fmd (t_{(y_B, my_F)}^{\star}E\otimes y^{\chi})=t_{\chi y_F}^{\star} \fmd (t_{(y_B, my_F)}^{\star}E)\otimes y_B^{\chi}\\ &=&t_{\chi y_F}^{\star} (t_{y_B}^{\star}\fmd (E)\otimes y_F^{-m})\otimes y_B^{\chi}\\ &=&t_{y_B+\chi y_F}^{\star} (I_{Z}^{\vee}(-r\sigma +\chi f )) \otimes \left(y_F^{-m}\boxtimes y_B^{\chi}\right)\\ &=&I_{{Z^-}}^{\vee}(-r\sigma + \chi f) \otimes y_F^{-d_v}\end{eqnarray*} where $Z^{-}=t_{y_B+\chi y_F}^{\star} Z,$ so that the addition in the fibers is $a_F(Z^-)=-\chi d_v y_F$. We therefore obtain the birational isomorphism $$\mathfrak M_v^{-}\dasharrow \mathfrak X_v^{-},\,\,\, \Phi_v^{-}(E, y)\mapsto (Z^-, -d_vy_F)$$ where $$\mathfrak X_v^{-}=\{(Z^{-}, z_F): a_F(Z^{-})=\chi z_F\}\hookrightarrow X^{[d_v]}\times F.$$ For further use, we record the following identities $$\fmd (V)=I_{Z^-}^{\vee}(-r\sigma+\chi f)\otimes z_F,\,\, \fmd(V^{\vee})=I_{\widetilde {Z^-}}(r\sigma-\chi f)\otimes z_F[-1].$$

\section {Rank-coprime arbitrary fiber degree} 
\label{arbitraryfiber}
Theorem \ref{t3} was proved in the previous Section for fiber degree $1$ in two ways: via an explicit analysis of O'Grady's description of the moduli space, and via Fourier-Mukai methods. In this section, we use Fourier-Mukai to prove Theorem \ref{t3} for arbitrary fiber degree. The argument builds on results of Bridgeland \cite {B}. At the end of the section, we briefly consider the case of a surface $X=C\times F$ with $F$ elliptic, but $C$ of higher genus. 

\subsection{Fourier-Mukai transforms in the general coprime setting} We write $$d=c_1(v)\cdot f$$  for the fiber degree, which we assume to be coprime to the rank $r$. Thus $c_1(v)=d\sigma+mf$. Pick integers $a$ and $b$ such that $$ad+br=1,$$ with $0<a<r$. The following lemma gives the kernel of the Fourier-Mukai transform we will use:

\begin {lemma} \label{l2}There exists a vector bundle $$\mathcal U\to F\times F$$ with the following properties:\begin {itemize}
\item [(i)] the restriction of $\mathcal U$ to $F\times \{y\}$ is stable of rank $a$ and degree $b$; 
\item [(ii)] the restriction of $\mathcal U$ to $\{x\}\times F$ is stable of rank $a$ and degree $r$.
\end {itemize}
Furthermore, $$c_1(\mathcal U)=b[o_F\times F] + r[F\times o_F] +c_1(\mathcal P_F).$$ 
\end {lemma}

\proof This result is known, see for instance \cite {B}, and can be explained in several ways. We consider the moduli space $M_F(a, b)$ of bundles of rank $a$ and degree $b$ over $F$. By the classic result of Atiyah \cite{a}, we have $$M_F(a, b)\cong F.$$ We let $$\mathcal U\to F\times F$$ denote the universal bundle. Therefore, for all $y\in F$, we have that $\mathcal U|_{E\times y}$ has rank $a$ and degree $b$; in fact the determinant equals $\mathcal O_F((b-1)\cdot o_F+y)$. The bundle $\mathcal U$ is not unique and we can normalize it in several ways. Indeed, for any matrix $$A=\begin {bmatrix} \lambda & a\\ \mu & b\end {bmatrix}\in SL_2(\mathbb Z),$$ we may assume that $\mathcal U|_{x\times F}$ has type $(a, \lambda)$. In fact, we can regard the first factor $F$ as the moduli space of bundles of rank $a$ and degree $\lambda$ over the second factor $F$, cf. \cite {B}. We may pick the pair $\lambda=r$ and $\mu=-d$. What we showed above allows us to conclude $c_1(\mathcal U)=b[o_F\times F] + r[F\times o_F] +c_1(\mathcal P_F)$.  Replacing $\mathcal U$ by a suitable twist, we may in fact achieve $$\det\, \mathcal U=\mathcal O(b[o_F\times F] + r[F\times o_F] )\otimes \mathcal P_F.$$ The lemma is proved.
 \qed

\begin {lemma} \label{l3}The sheaf $\mathcal U\to F\times F$ is semihomogeneous. More precisely, 
\begin{equation}
\label{translate}
t_{(x, y)}^{\star} \mathcal U=\mathcal U \otimes \pi_1^{\star} (x^{-\frac{b}{a}}y^{\frac{1}{a}})\otimes \pi_2^{\star}(x^{\frac{1}{a}}y^{-\frac{r}{a}}).
\end{equation}
In particular, $$\text{ch }\mathcal U=a \exp \left(\frac {c_1(\mathcal U)}{a}\right)\implies \chi(\mathcal U)=-d.$$
\end {lemma}
\proof By symmetry it suffices to argue that $$t_{(x, 0)}^{\star} \mathcal U=\mathcal U\otimes \pi_1^{\star} x^{-\frac{b}{a}}\otimes \pi_2^{\star} x^{\frac{1}{a}}.$$ We note that the choice of roots  for the line bundles on the right hand side is not relevant. The restrictions of both sides to $F\times \{y\}$ agree: $\mathcal U_y=\mathcal U|_{F\times \{y\}}$ is stable on $F$, and thus semihomogeneous, satisfying $$t^{\star}_x \,\mathcal U_y=\mathcal U_y\otimes x^{-\frac{b}{a}}.$$ We check agreement over $o_F\times F$. This is the statement that $$\mathcal U|_{x\times F}=\mathcal U|_{o\times F}\otimes x^{\frac{1}{a}}$$ which holds by comparing ranks and determinants. Finally, agreement over $F\times F$ follows from the generalized see-saw Lemma $2.5$ of Ramanan \cite {R}. 

The second part of the lemma concerning the numerical invariants of $\mathcal U$ follows from general facts about semihomogeneous bundles \cite {muk}. \qed
\vskip.1in

\begin {remark} A family of semihomogeneous bundles with fixed numerical invariants were constructed in arbitrary dimension in \cite {O}, and played a role in the decomposition of the Verlinde bundles. The dimension $1$ case specializes to the bundle $\mathcal U$ considered here. 
\end {remark}
Letting $\pi_{F\times F}: F\times F\times B\to F \times F$ be the projection, we consider the Fourier-Mukai transform 
$$\fmd:\mathbf D(X)\to \mathbf D(X)$$  with kernel $$\pi^{\star}_{F\times F}\, \mathcal U\to F\times F\times B\cong X\times_{B}X.$$ It follows from \cite {B} that the kernel $\mathcal U$ is strongly simple over each factor, hence the Fourier-Mukai transform $\fmd$ is an equivalence, with inverse having kernel $\pi_{F\times F}^{\star}\, \mathcal U^{\vee}[1]$.
\begin {lemmaa}
$$\fmd (E\otimes y)=t_{a y_F}^{\star} \fmd (E)\otimes (y_B \boxtimes y_F^{r}),$$
$$\fmd (t_x^{\star}E)=t_{x_B+bx_F}^{\star} \fmd (E)\otimes x_F^{-d}.$$
\end {lemmaa}
\proof For both formulas, it is enough to consider the case of split sheaves $$E=G\boxtimes H,$$ where $G, H$ are sheaves over $B$ and $F$ \footnote{Indeed, both left and right hand sides of the equalities claimed by the Lemma are Fourier-Mukai equivalences (note that translations and tensorization are particular examples of Fourier-Mukai transforms, and composition of Fourier-Mukai transforms is one as well). Taking inverses, it suffices to prove that if a Fourier-Mukai equivalence is the identity over split sheaves of the type $G\boxtimes H$, then it is always the identity. But this is clear, as one proves that the kernel of such a transform is the structure sheaf of the diagonal.}.  Clearly, $$\fmd (E)=G \boxtimes \rp \pi_{2\star}(\pi_1^{\star} H\otimes \mathcal U),$$ while$$\fmd(E\otimes y)=\left(G \otimes y_B\right) \boxtimes \rp \pi_{2\star}(\pi_1^{\star}(H\otimes y_F)\otimes \mathcal U).$$ For the first formula, it suffices now to explain that \begin {equation}\label{e3}\rp \pi_{2\star}(\pi_1^{\star}(H\otimes y_F)\otimes \mathcal U)=t_{ay_F}^{\star} \rp \pi_{2\star}(\pi_1^{\star}H\otimes \mathcal U)\otimes y_F^{{r}}.\end{equation} This is however clear since from \eqref{translate}, 
$$\mathcal U\otimes \pi_1^{\star}y_F=t_{(0, ay_F)}^{\star} \mathcal U\otimes \pi_2^{\star} y_F^{{r}}.$$ We now explain the second formula by assuming as before that $E$ splits. We calculate $$\fmd (t_x^{\star}E)=\fmd(t_{x_B}^{\star} G \boxtimes
 t_{x_F}^{\star}H)= t_{x_B}^{\star} G \boxtimes \rp \pi_{2\star}(\pi_1^{\star}t_{x_F}^{\star} H\otimes \mathcal U).$$ It remains to argue that $$\rp \pi_{2\star}(\pi_1^{\star}t_{x_F}^{\star} H\otimes \mathcal U)=t_{bx_F}^{\star}\rp \pi_{2\star}(\pi_1^{\star}H\otimes \mathcal U)\otimes x_F^{-d}.$$ Indeed, \eqref{translate} gives $$t_{(x_F, 0)}^{\star}\mathcal U=\mathcal U\otimes \pi_1^{\star}x_F^{-\frac{b}{a}}\boxtimes \pi_2^{\star} x_F^{\frac{1}{a}},$$ so we calculate \begin{eqnarray*} \rp \pi_{2\star} (\pi_1^{\star}t_{x_F}^{\star} H\otimes \mathcal U )&=&\rp \pi_{2\star} \left (t_{(x_F, 0)}^{\star}(\pi_1^{\star}H\otimes \mathcal U)\otimes \pi_1^{\star} x_F^{\frac{b}{a}} \otimes \pi_2^{\star}x_F^{-\frac{1}{a}}\right )\\ &=&\rp \pi_{2\star}\left(t_{(x_F, 0)}^{\star}\left (\pi_1^{\star}(H\otimes x_F^{\frac{b}{a}})\otimes \mathcal U\right ) \right)
\otimes x_F^{-\frac{1}{a}}\\ &=& \rp \pi_{2\star}\left(\pi_1^{\star}(H\otimes x_F^{\frac{b}{a}})\otimes\mathcal U\right)\otimes x_F^{-\frac{1}{a}}\\ &=&t_{bx_F}^{\star}\rp \pi_{2\star}(\pi_1^{\star} (H)\otimes \mathcal U)\otimes x_F^{\frac{br}{a}} x_F^{\frac{-1}{a}}\text { (using } \eqref{e3})\\ &=&t_{bx_F}^{\star}\rp \pi_{2\star}(\pi_1^{\star} (H)\otimes \mathcal U)\otimes x_F^{-d}.\end{eqnarray*}\qed

\begin {propositiona} 
\label{pa}
For a generic sheaf $V$ of fixed determinant and determinant of Fourier-Mukai $$\det V=\mathcal O(d\sigma + m f),\,\, \det \widehat V=\mathcal O(-m\sigma - df)$$ we have $$\fmd(V)=I_Z^{\vee}\otimes \mathcal O_X((a\chi+bm)f),$$ for a subscheme $Z$ of length $d_v$ with $a(Z)=0$. 
\end{propositiona}

\proof In the proof of this Proposition, it will be important to distinguish between the two copies of $X$ which are the source and the target of the Fourier-Mukai transform, because of the asymmetry present in the bundle $\mathcal U$. We will write $X_1$ and $X_2$ for two copies of $X$, respectively. 

By Grothendieck duality, to prove that $$\fmd(V)=I_Z^{\vee}\otimes \mathcal O_X((a\chi+bm)f),$$ it suffices to show that $$\Psi_{X_1\to X_2}^{\mathcal U^{\vee}} (V^{\vee})=I_Z\otimes \mathcal O(-(a\chi+bm) f)[-1],$$ where $\Psi$ is the Fourier-Mukai transform with kernel $\mathcal U^{\vee}$. Using \cite {B}, Lemma $6.4$, we know that $V^{\vee}$ is $WIT_1$ with respect to $\Psi$, since the restriction to the general fiber is stable, of slope $-d/r<b/a$. Note that $\Psi(V^{\vee})[1]$ has rank $$-\chi(V^{\vee}|_{F\times y}\otimes \mathcal U^{\vee}|_{F\times y})=ad+br=1.$$ Section $7$ of \cite {B}, or Sections $3$ and $5$ in \cite {BH}, show that for generic $V$, the $\Psi$-transform is torsion-free, hence it must be of the form $I_Z\otimes L[-1]$ for some line bundle $L$. Bridgeland's argument moreover shows that the subscheme $Z$ has length $d_v$. 

The fiber degree of $\Psi(V^{\vee})$ equals $$c_1(\rp \pi_{2\star}(\pi_1^{\star}V|^{\vee}_{\mathsf f}\otimes \mathcal U^{\vee}))= \pi_{2\star}(\pi_1^{\star}(r-d\omega) (a-c_1(\mathcal U)+\text{ch}_2(\mathcal U)))_{(2)}=0,$$ where Lemmas \ref{l2} and \ref{l3} are used to express the numerical invariants of $\mathcal U$. In fact more is true. Since the restriction of $V$ to a generic fiber is stable, it must equal the Atiyah bundle $\mathsf E_{r, d}$. This implies that the restriction of $L$ to a generic fiber must coincide with $\Psi (\mathsf E_{r, d})[1]$ which is trivial. Therefore, $L$ must be a sheaf of the form $\pi_B^{\star} M^{\vee}$, where $M$ is a degree $-\beta$ line bundle over $B$. We prove $$\beta=-a\chi -bm.$$

To this end, we calculate
the Euler characteristic of $\Psi(V^{\vee}(\sigma))$ as \begin{eqnarray*}\chi(\Psi(V^{\vee})(\sigma))&=&\chi(F\times F\times B, \pi_{13}^{\star}V^{\vee}\otimes \pi_{12}^{\star} \mathcal U^{\vee}\otimes \pi_{23}^{\star} \mathcal O(\sigma))\\&=& \int_{F\times F\times B} \pi_{13}^{\star}(r-(d\sigma+mf)+\chi \omega)\cdot \pi_{12}^{\star}(a-c_1(\mathcal U)+\text{ch}_2(\mathcal U))\cdot \pi_{23}^{\star}(1+\sigma)\\&=&bm+a\chi-\chi r +md.\end{eqnarray*} On the other hand, $$\chi(I_Z\otimes L(\sigma))=\beta-d_v=\beta - (dm-r\chi)$$ hence $\beta=-a\chi-bm.$

When the determinant and determinant of Fourier-Mukai of $V$ are fixed, we show that $a(Z)=0$. First, we analyze the requirement the determinant be fixed. The inverse of $\fmd$ is given by $\Phi_{X_2\to X_1}^{\mathcal U^{\vee}[1]},$ the Fourier-Mukai whose kernel is $\mathcal U^{\vee}[1]$, considered as a transform from $X_2\to X_1$. Hence, $$V=\Phi_{X_2\to X_1}^{\mathcal U^{\vee}[1]} (L^{\vee}\otimes I_Z^{\vee}), \,\,\, L=\pi_B^{\star} M^{\vee}$$ has fixed determinant $\mathcal O(d \sigma+mf)$. In order to make the computations more explicit, we write $$M=\mathcal O_B(-(\beta+1) [o_B]+[\mu]),$$ for some $\mu\in B$. We have \begin{eqnarray}\nonumber\det V^{\vee} &=&\det \Phi^{\mathcal U^{\vee}}(L^{\vee})\bigotimes_{z\in Z}\det \rp \pi_{13!}(\pi_{12}^{\star}\mathcal U^{\vee}\otimes \pi_{23}^{\star}\mathcal O_z^{\vee})^{\vee}\\\nonumber&=&\det \Phi^{\mathcal U^{\vee}}(L^{\vee})\otimes \bigotimes \det \left(\mathcal U^{\vee}|_{F\times z_F}\boxtimes \mathcal O_{z_B}[2]\right)^{\vee}\\ &=&\nonumber \det \Phi^{\mathcal U^{\vee}}(L^{\vee})\otimes_{z\in Z} (\mathcal O_F\boxtimes \mathcal O_B(-a [z_B]))\\&=&\nonumber\det \Phi^{\mathcal U^{\vee}}(L^{\vee})\otimes \left(\mathcal O_F\boxtimes  \mathcal O_B(-(ad_v-1) [o_B]-[a\cdot a_B(Z)])\right)\\&=&\label{fl}\mathcal O(-d\sigma-mf)\otimes \pi^{\star} \mathcal O_B(-[a\cdot a_B(Z)+r\mu]+[o_B]).\end{eqnarray} This gives $$a\cdot a_B(Z)+r\mu=o_B.$$ In equation \eqref{fl}, we used the calculation $$\det \Phi^{\mathcal U^{\vee}}(L^{\vee})=\det \rp p_{13!}(p_{12}^{\star}\mathcal U^{\vee}\otimes p_3^{\star} M)=\det \left(\rp p_{1!}\mathcal U^{\vee}\boxtimes M\right)=\det \rp p_{1!} (\mathcal U^{\vee}) \,\boxtimes M^{-r}$$ $$=\mathcal O_F(-do_F)\boxtimes \mathcal O_B((r\beta+1)[o_B]-[r\mu]),$$ where by Lemma $3$, the pushforward of $\mathcal U^{\vee}$ has rank $-r$ and degree $-d$. Equation \eqref{fl} also makes use of the identity $$r\beta-ad_v=-m.$$

We analyze the requirement that the determinant of the Fourier-Mukai be fixed. We know that $$\fm (V)=\fm \circ \Phi_{X_2\to X_1}^{\mathcal U^{\vee}[1]} (L^{\vee}\otimes I_Z^{\vee}).$$ This composition can be re-expressed as a Fourier-Mukai whose kernel equals the convolution of the following two kernels: $\widetilde{\mathcal U}^{\vee}[1]$ for $\Phi$, and $$\mathcal P=\mathcal P_F\times \mathcal P_B$$ for $\fm$. The tilde sign indicates that the kernel $\mathcal U^{\vee}[1]$ is considered in the opposite direction for $\Phi$ than it is for $\Psi$. This is the same as applying to  $\mathcal U\to F\times F$the involution that exchanges the factors. The new kernel can be expressed as $$\rp p_{13\star}(p_{12}^{\star}\widetilde{\mathcal U}^{\vee}[1]\otimes p_{23}^{\star} \mathcal P)=\rp \pi_{13\star}(\pi_{12}^{\star}\widetilde {\mathcal U}^{\vee}\otimes \pi_{23}^{\star}\mathcal P_F)[1]\boxtimes \mathcal P_B=\mathcal V\boxtimes \mathcal P_B$$ where $$\mathcal V\to F\times F,\,\, \mathcal V=\rp \pi_{13\star}(\pi_{12}^{\star}\widetilde{\mathcal U}^{\vee}\otimes \pi_{23}^{\star}\mathcal P_F)[1]$$ is the fiberwise Fourier-Mukai image of $\widetilde{\mathcal U}^{\vee}$ up to a shift. 
The complex $\mathcal V$ has rank $b$, and a Riemann-Roch calculation shows that $$c_1(\mathcal V)=d[o_F\times F]+a[F\times o_F]+c_1(\mathcal P_F).$$ The pushforward $\rp (p_2^{F})_{\star}( \mathcal V)$ has rank $d$ and determinant $-r[o_F]$. This will be used below. Fiberwise, note that $$\det \mathcal V^{\vee}|_{z_F\times F}=\det \rp p_{2!} (p_1^{\star}\widetilde {\mathcal U}|_{z_F\times F}^{\vee}\otimes \mathcal P_F)=-[z_F]-(a-1)[o_F].$$ 

Now, we calculate \begin{eqnarray*}\det \fm (V)&=&\det \rp p_{2\star}(\mathcal V \boxtimes \mathcal P_B\otimes p_1^{\star}(L^{\vee}\otimes I_Z^{\vee}))\\&=&\det \rp p_{2\star}(\mathcal V\boxtimes\mathcal P_B\otimes p_1^{\star}L^{\vee})\otimes \det \rp p_{2\star} (\mathcal V\boxtimes \mathcal P_B\otimes p_1^{\star}\mathcal O_Z[2])^{\vee}.\end{eqnarray*}
The first determinant is constant $$\det \rp p_{2\star}(\mathcal V \boxtimes\mathcal P_B \otimes p_1^{\star}L^{\vee})=\det \left(\rp (p^F_{2})_{\star} (\mathcal V) \boxtimes \rp (p_2^B)_{\star}(\mathcal P_B\otimes (p_{1}^B)^{\star} M \right)) $$ $$=\left(\det \rp (p_2^{F})_{\star}( \mathcal V)\right)^{-\beta}\boxtimes \mathcal O_B(-[-\mu])^{d}=\mathcal O_F(r\beta [o_F])\boxtimes \mathcal O_B(-d[-\mu]).$$
The second determinant needs to be fixed, and it equals $$\bigotimes_{z\in Z}\det \mathcal V^{\vee}|_{z_F\times F} \boxtimes \mathcal O_B(-[z_B]+[o_B])^{b}=\bigotimes_{z\in Z}\mathcal O_F(-[z_F]-(a-1)[o_F])\boxtimes \mathcal O_B(-[z_B]+[o_B])^{b}$$ $$=\mathcal O_F(-[a_F(Z)]-(ad_v-1)[o_F])\boxtimes \mathcal O_B(-[b\cdot a_B(Z)]+[o_B]).$$ Therefore $$\det \fm (V)=\mathcal O(-mf-d\sigma)\otimes \mathcal O_F(-[a_F(Z)]+[o_F])\boxtimes \mathcal O_B(-[b\cdot a_B(Z)-d\mu]+[o_B]).$$ Since $$\det \fm(V)=\mathcal O(-mf-d\sigma),$$ this immediately yields $$a_F(Z)=o_F \text { and } b\cdot a_B(Z)-d\mu=o_B.$$ 

Combining these two equations with $$a\cdot a_B(Z)+r\mu=o_B$$ shown aprove, and the fact that $ad+br=1$, we obtain $a_B(Z)=0$ and $\mu=0.$ Thus $$a(Z)=0,\,\, L=\mathcal O((a\chi+bm)f).$$ 
This completes the proof. 
\qed
\vskip.1in

\noindent {\it Proof of Theorem \ref{t3}.} In the course of the proof, we showed that for a generic sheaf $V$ of rank $r$ and determinant $$\det V=\mathcal O(d\sigma + m f),$$ the Fourier-Mukai transform takes the form $$\fmd(V)=I_Z^{\vee}((a\chi+bm)f)\otimes \pi_{B}^{\star}\mu,$$ for some $\mu\in B\cong \widehat B$ such that $$a\cdot a_B(Z)+r\mu=o_B.$$ The assignment $$V\to (Z, ba_B(Z)-d\mu)$$ gives the birational isomorphism $$\mathfrak M_v^+\dasharrow \mathfrak X_v^+,$$ claimed by Theorem \ref{t3}. 

\begin {remark} In a similar fashion we could find the birational type of the moduli space $\mathfrak M_v^{-}.$ However, it is not necessary to repeat the argument above to deal with this new case. We could instead make use of the map
$$\Phi_v^-:K_v\times X\to \mathfrak M_v^-, \,\, (E, y)\mapsto t^{\star}_{(dy_B, my_F)} E\otimes y^{\chi}$$ and use Lemma 1A to calculate \begin{eqnarray*}\fmd (\Phi_v^{-}(E, y))&=&\fmd (t_{(dy_B, my_F)}^{\star}E\otimes y^{\chi})=t_{a\chi y_F}^{\star} \fmd (t_{(dy_B, my_F)}^{\star}E)\otimes (y_B^{\chi}\boxtimes y_F^{r\chi})\\ &=&t_{a\chi y_F}^{\star} (t_{dy_B+bmy_F}^{\star}\fmd (E)\otimes y_F^{-dm})\otimes (y_B^{\chi}\boxtimes y_F^{r\chi})\\ &=&t_{dy_B+(a\chi+bm) y_F}^{\star} (I_{Z}^{\vee}((a\chi+bm) f )) \otimes \left(y_B^{\chi}\boxtimes y_F^{-d_v}\right)\\ &=&I_{{Z^-}}^{\vee}((a\chi+bm)f) \otimes \left(y_B^{\chi - (a\chi+bm)d}\boxtimes y_F^{-d_v}\right)\\&=&I_{{Z^-}}^{\vee}((a\chi+bm)f) \otimes \left(y_B^{-bd_v}\boxtimes y_F^{-d_v}\right)\end{eqnarray*} where $Z^{-}=t_{dy_B+(a\chi+bm) y_F}^{\star} Z.$ Thus, the assignment $$V\mapsto (Z^-, -d_vy)$$ gives a birational isomorphism $$\mathfrak M_v^-\dasharrow\mathfrak X_v^-=\{(Z^-, z): a(Z^-)=f(z)\}\hookrightarrow X^{[d_v]}\times X$$ where the isogeny $f:X\to X$ is given by $$f(z)=(dz_B, (a\chi+bm)z_F).$$ In the case of fiber degree $1$, this specializes to the subvariety $$\mathfrak X_v^-=\{(Z^-, z_F): a_F(Z^-)=\chi z_F\}\hookrightarrow X^{[d_v]}\times F$$ of Section \ref{s2}.
\end {remark}

\subsection {Higher genus} \label{g2} Assume now $(C, o)$ is a pointed smooth curve of genus $g\geq 1$, and $F$ is still an elliptic curve. We set $\gbar=g-1$. Consider the product surface $$X=C\times F\to C.$$ Let $\mv$ be the moduli space of sheaves over $X$ of rank $r$ and determinant $\mathcal O(\sigma +m f)$, where $\sigma$ is the zero section and $f$ denotes the fiber over $o$. We describe the birational type of $\mv$, in codimension $1$ for $r\neq 2$, using the Fourier-Mukai transform with kernel the Poincar\'e bundle $$\pi_{F\times F}^{\star} \mathcal P_F\to F\times F\times C.$$ The proof is entirely similar to that of Proposition 1A, so we content here to only record the result. 

Let $$a_C:X^{[d_v]}\to C^{[d_v]}$$ be the map induced by the projection $X\to C$. Thus, each scheme $Z$ of length $d_v$ in $X$ yields a divisor $a_C(Z)$ of degree $d_v$ over the curve $C$. The line bundle $$\mathcal M_Z=\mathcal O_C(a_C(Z)-d_v\cdot o)$$ has degree $0$, and therefore admits roots of order $r$. We define $$\mathfrak X_v^+=\{(Z, c): c^r=\mathcal M_Z\}\hookrightarrow X^{[d_v]}\times \text{Pic}^0(C).$$ Then, for $V\in \mv$, we have $$\fmd (V)=I_Z^{\vee}(-r\sigma+(\chi+\gbar) f)\otimes c^{-1},$$ establishing the birational isomorphism $$\mathfrak M_v^+\dasharrow \mathfrak X_v^+.$$ The same statement holds in any fiber degree coprime to the rank, but we will not detail this fact. 


\section{The strange duality isomorphism} 
We now proceed to prove Theorems \ref{t1}, \ref{t2} and $2$A stated in the introduction. Throughout this section, we place ourselves in the context when the fiber degree is $1$. 
	
\subsection {Reformulation}

Let $X=B\times F$ be a product abelian surface. As a consequence of Section \ref{s2}, under the birational map $$\mv \times \mw \dasharrow \mathfrak X^{+}_v \times \mathfrak X^{+}_w$$ induced by the relative Fourier-Mukai transform, the standard theta divisor 
$$\Theta^+_{vw} = \{(V, \, W) \, \, \text{with} \, \, h^1 (V \otimes W) \neq 0 \} \subset \mv \times \mw$$
is identified with a divisor $$\Theta^{+}\subset \mathfrak X_v^{+}\times \mathfrak X_w^{+}.$$
Note that for sheaves $(V, W)\in \mathfrak M_v^+\times \mathfrak M_w^+$ corresponding to pairs $$(Z^+, z_B)\in \mathfrak X_v^+ \text{ and }(T^+, t_B)\in \mathfrak X_w^+,$$ we have \begin{eqnarray*}H^1(V\otimes W)&=&\text{Ext}^1(W^{\vee}, V)=\text{Ext}^1(\fmd(W^{\vee}), \fmd (V))\\&=&\text{Ext}^1 ( I_{\widetilde {T^{+}}} (s\sigma - \chi' f)[-1]\otimes t_B, I_{Z^+}^{\vee}\otimes \mathcal O(-r\sigma +\chi f) \otimes z_B^{-1})\\ &=&\text{Ext}^1(I_{Z^+}^{\vee}\otimes \mathcal O(-r\sigma +\chi f) \otimes z_B^{-1}, 
I_{\widetilde {T^{+}}} (s\sigma - \chi' f)[-1]\otimes t_B)^{\vee}\\ &=&H^0(I_{Z^+}\otimes I_{\widetilde{T^+}}\otimes z_B \otimes t_B\otimes \mathcal O((r+s)\sigma -(\chi+\chi')f))^{\vee}.\end{eqnarray*} (The notation above has the obvious meaning: $r, s$ are the ranks of $v$ and $w$, while $\chi, \chi'$ are their Euler characteristics.) Thus, the theta divisor $\Theta^+_{vw}$ in the product $\mv\times \mw$ corresponds to the divisor $$\Theta^{+}=\{(Z^+, z_B, T^+, t_B): h^0(I_{Z^+}\otimes I_{\widetilde {T^+}}\otimes z_B \otimes t_B \otimes L)\neq 0\}$$ in the product $\mathfrak X_v^+\times \mathfrak X_w^+$.  
Here, we set \begin{equation}\label{expl}L=\mathcal O((r+s)\sigma -(\chi+\chi')f)\, \, \, \text{on} \, \, \, X.\end{equation} In consequence, strange duality is demonstrated if we show that the divisor $\Theta^{+}$ induces an isomorphism $$\mathsf D^+: H^0 \left (\mathfrak X^{+}_v, \Theta_w \right )^{\vee} \longrightarrow H^0 \left (\mathfrak X^{+}_w, \Theta_v \right ).$$ 

\subsection {Theta bundles} Since for $r\neq 2$ the birational isomorphism $\mv  \dasharrow \mathfrak X^{+}_v$ is defined away from codimension 2, we are interested in an explicit description of the determinant line bundle  $$\Theta_w\to \mathfrak X_v^{+}.$$ For a line bundle $L$ on $X$, we standardly let $$L^{[d_v]} = \det Rp_{\star}\left (  {\mathcal O}_{\mathcal Z} \otimes q^{\star} L \right ) \, \, \, \text{on} \, \, \, X^{[d_v]}.$$ We also note the natural projections $$c_v:\mathfrak X_v^+\to X^{[d_v]}, \,\, (Z^+, z_B)\mapsto Z^+$$ and $$\pi_2:\mathfrak X_v^+\to B,\,\, (Z^+, z_B)\to z_B.$$ The theta bundle is calculated by the following result:

\begin {proposition}
\label{basicid} We have 
$$\Theta_w=c_v^{\star} L^{[d_v]} \otimes \pi_2^{\star} {\mathcal O}_B ((s-r)o_B),$$ with $L$ given by \eqref{expl}. 
\end {proposition}

\proof We give one proof of the proposition here; another one is essentially contained in Section 4.4. We begin by noting the degree $d_v^4$ \'{e}tale morphism 
$$q_v: K^{[d_v]} \times X = K^{[d_v]}\times B \times F \to \mathfrak X_v^{+},\,\,\, (Z, x)\mapsto (t_{rx_B+x_F}^{\star} Z, \, -d_v x_B)$$
It is related to the standard less twisted map 
$$\mu_v:K^{[d_v]}\times X\to X^{[d_v]}, \, \, \, \, \mu_v(Z, \, \, x)=t_{x}^{\star} Z$$
via the commutative diagram
\begin{center}
$\xymatrix{ {K^{[d_v]} \times B \times F}  \ar[d]^{(1,\, r,\, 1)}  \ar[r]^{q_v} & \mathfrak X^{+}_v \ar[d]^{c_v} \\ K^{[d_v]} \times B \times F\ar[r]^{\mu_v} & {X^{[d_v]}}.}$
\end{center}
Using the diagram, we calculate
$$q_v^{\star}\, c_v^{\star} L^{[d_v]} =(1, r, 1)^{\star} \mu_v^{\star} L^{[d_v]} =(1, r, 1)^{\star} \left (L^{[d_v]}\boxtimes L^{d_v}\right )=L^{[d_v]}\boxtimes \left ( (r, 1)^{\star} L\right ) ^{d_v},$$
and further,
$$(r, 1)^{\star} L= (r,1)^{\star} {\mathcal O} \left ( (r+s ) \sigma - (\chi + \chi') f \right ) =  {\mathcal O} \left ( (r+s ) \sigma - r^2 (\chi + \chi') f \right ).$$
If $$\pi_B: K^{[d_v]} \times B \times F \to B$$ is the projection to $B$, from the definitions we also have $\pi_2\circ q_v=-d_v$ hence 
$$q_v^{\star} \, \pi_2^{\star}  {\mathcal O}_B ((s-r)o_B) = \pi_B^{\star}  {\mathcal O}_B \left ( (s-r) d_v^2o_B \right ).$$ 

Putting the previous three equations together we find
\begin{equation} \label{pb1} q_v^{\star}\left (  c_v^{\star} L^{[d_v]} \otimes \pi_2^{\star} {\mathcal O}_B ((s-r)o_B) \right )  =  L^{[d_v]} \boxtimes {\mathcal O} \left ((r+s) \sigma +( (s-r ) d_v - r^2 (\chi + \chi') ) f \right )^{ \otimes d_v}.\end{equation}
On the other hand, Proposition 2 of  \cite {abelian} (written in holomorphic $K$-theory, as in \cite {O}) gives
\begin{equation} \label{pb2} q_v^{\star}\,  \Theta_w = L^{[d_v]} \boxtimes {\mathcal O} \left ( (r+s) \sigma + (rn + sm) f \right )^{\otimes d_v}.\end{equation}
The two pullbacks \eqref{pb1} and \eqref{pb2} are seen equal on $K^{[d_v]} \times X$, as one shows that $$rn+sm=(s-r)d_v-r^2(\chi+\chi').$$ This uses the numerical identity $$m+n = - r\chi' - s\chi$$ which expresses strange duality orthogonality $\chi (v\cdot w)=0,$ also remembering that $d_v=m- r\chi.$

Now, consider $$Q=\Theta_w\otimes \left(c_v^{\star} L^{[d_v]}\right)^{\vee}\otimes \pi_B^{\star} \mathcal O_B((r-s)o_B).$$ We showed that $q_v^{\star} Q$ is trivial. Note the morphism $$p_v:\mathfrak X_v^+\to X,\,\, (Z, x_B)\to (a_F(Z), x_B)$$ with fibers isomorphic to $K^{[d_v]}$. In fact, each fiber $$\iota: p_v^{-1}(x)\hookrightarrow \mathfrak X_v^+$$ factors through the morphism $$q_v:K^{[d_v]}\times X\to \mathfrak X_v^+.$$ Indeed, $\iota=q_v\circ j,$ where $j: p_v^{-1}(x)\to K^{[d_v]}\times X$ is the map $$j(Z)=(t_{-ry_B-y_F}^{\star} Z, y),$$ for any choice of $y\in X$ such that $d_v y=-x$. Therefore, the above argument implies that the restriction of $Q$ to each fiber $p_v^{-1}(x)$ is trivial. Hence, $$Q=p_v^{\star} N,$$ for some line bundle $N$ over $X$. 

We now argue that $N$ is trivial, by constructing a suitable test family. Consider $Z_0$ a subscheme of length $d_v-1$ supported at $0$, and define $$\alpha:X \dasharrow
\mathfrak X_v^{+},\,\,\, x\mapsto (Z_0+(rx_B, x_F), x_B).$$ The map $\alpha$ is defined away from the $r^2$ points in $B[r]\times o_F.$ It suffices to show $N$ is trivial along this open set. Pulling back the equality
$$p_v^{\star}N=\Theta_w\otimes \left(c_v^{\star} L^{[d_v]}\right)^{\vee}\otimes \pi_B^{\star} \mathcal O_B((r-s)o_B)$$ under $\alpha$, and noting $\alpha\circ p_v=1$, it suffices to prove that \begin{equation}\label{test}\alpha^{\star} \Theta_w=\overline c_v^{\star} L^{[d_v]}\otimes \pi_B^{\star}\mathcal O_B((s-r)o_B)\end{equation} where $$\overline c_v=c_v\circ \alpha: X\dasharrow X^{[d_v]} \text{ is given by } \overline c_v(x)=Z_0+(rx_B, x_F).$$ 

We calculate the right hand side. The universal family $$\mathcal Z\subset X^{[d_v]}\times X$$ becomes under pullback by $\overline c_v$ the family $Z_0+\Delta_r\cong \Delta_r$, where $\Delta_r\subset X\times X$ is the $r$-fold diagonal $$\Delta_r=(rx_B, x_F, x_B, x_F).$$ Then,  $$\overline c_v^{\star}L^{[d_v]}=\overline c_v^{\star} \det \rp p_! (\mathcal O_{\mathcal Z}\otimes q^{\star} L)=\det \mathbf Rp_!( \mathcal O_{\Delta_r}\otimes q^{\star}L)=(r_B, 1_F)^{\star}L$$ hence $$\overline c_v^{\star}L^{[d_v]}=\mathcal O((r+s)\sigma-r^2(\chi+\chi')f).$$

For the left hand side, fix $(T, c)$ in $\mathfrak X_w^{+}$, where $c\in B\cong \widehat B$, so that $$a_B(T)=s c.$$ We have $$\Theta_w=\det \rp p_!(I_\mathcal Z\otimes \mathcal P_B\otimes q^{\star}(L\otimes c\otimes I_{\widetilde T}))^{\vee},$$ where $\mathcal P_B$ is the Poincar\'e bundle over $B\times B$. Write $$M=L\otimes c \otimes I_{\widetilde T}$$ so that $$M=M_B\boxtimes M_F-\sum_{t\in T} \mathcal O_{t_B\times -t_F}$$ for $$M_B=c \otimes \mathcal O_B(-(\chi+\chi')o_B), \,\, M_F=\mathcal O((r+s)o_F).$$ Since the universal family pulls back to $\alpha^{\star} \mathcal Z\cong \Delta_r$, we have that $$\alpha^{\star} I_{\mathcal Z}=\mathcal O-\mathcal O_{\Delta_r}.$$ Therefore, $$\alpha^{\star} \Theta_w=\det \rp p_!(\mathcal P_B\otimes q^{\star}(M_B\boxtimes M_F))^{\vee} \otimes  \det \rp p_{!} (\mathcal O_{\Delta_r}\otimes \mathcal P_B\otimes q^{\star}(M_B\boxtimes M_F))$$
$$\otimes_{t\in T} \left(\det \rp p_!(\mathcal P_B\otimes q^{\star}\mathcal O_{t_B\times -t_F})\otimes \det \rp p_!(\mathcal O_{\Delta_r}\otimes \mathcal P_B\otimes q^{\star} \mathcal O_{t_B\times -t_F})^{\vee}\right).$$ We calculate the first term $$\det \rp p_{!} (\mathcal P_B\otimes q^{\star} (M_B\boxtimes M_F))^{\vee}=\det (\rp p_{B!}(\mathcal P_B\otimes q_B^{\star} M_B )\boxtimes H^{\bullet}(M_F)\otimes \mathcal O_F)^{\vee}$$ $$=\det (\rp p_{B!}(\mathcal P_B \otimes q_B^{\star} M_B ))^{-(r+s)}\boxtimes \mathcal O_F= \left(c\otimes \mathcal O_B(-o_B)\right)^{-(r+s)}\boxtimes \mathcal O_F.$$ The second term becomes 
$$\det \rp p_! (\mathcal O_{\Delta_r}\otimes \mathcal P_B \otimes q^{\star} (M_B\boxtimes M_F))=\det (\rp p_{B!}(\mathcal O_{\Delta^B_r}\otimes \mathcal P_B\otimes q_B^{\star} M_B)\boxtimes M_F)$$ where $\Delta_r^B$ is the image of $$j:B\to \Delta_r^B, \,\, x_B\to (rx_B, x_B).$$ We calculate $$\rp p_{B!}(\mathcal O_{\Delta_r^{B}}\otimes (\mathcal P_B\otimes q_B^{\star} M_B))=j^{\star} \mathcal P_B\otimes r^{\star} M_B=\mathcal O_B(-2r o_B)\otimes r^{\star} M_B.$$ Therefore, the second term equals
$$c^{r}\otimes \mathcal O_B(-2ro_B-r^2(\chi+\chi')o_B)\boxtimes \mathcal O_F((r+s)o_F).$$ 
The third term now equals $$\det \rp p_!(\mathcal P_B \otimes q^{\star} \mathcal O_{t_B\times -t_F})=\mathcal P_{t_B}\boxtimes \mathcal O_F,$$ hence the tensor product over all $t\in T$ yields $$\mathcal P_{a_B(W)}\boxtimes \mathcal O_F=c^{s}\boxtimes \mathcal O_F.$$ Finally, the fourth term is easily seen to be trivial $$\det \rp p_! (\mathcal O_{\Delta_r}\otimes \mathcal P_B \otimes q^{\star} \mathcal O_{t_B\times -t_F})=\mathcal O_X.$$ Equation \eqref{test} follows putting all the terms together. This concludes the proof of Proposition \ref{basicid}.
\qed
 \vskip.1in
 
\subsubsection {Fixed determinant of Fourier-Mukai} The discussion for the moduli space of sheaves with fixed determinant of their Fourier-Mukai is entirely parallel. We identify the theta divisor $$\Theta^{-}_{vw}\subset \mathfrak M_v^{-}\times \mathfrak M_w^{-}$$ with the divisor $$\Theta^{-}=\{(Z^-, z_F, T^-, t_F): h^0(I_{Z^-}\otimes I_{\widetilde T^-}\otimes z_F^{-1}\otimes t_F \otimes L)\neq 0\}\hookrightarrow \mathfrak X_v^-\times \mathfrak X_w^-,$$ where as before $$L=\mathcal O(-(\chi + \chi')f+(r+s)\sigma).$$ For $r, s\geq 3$, the birational isomorphisms are defined in codimension $1$, and we have $$\Theta_w=(c_v^{-})^{\star} L^{[d_v]} \otimes \pi_{F}^{\star} \mathcal O_F((\chi - \chi')o_F),$$ where $$c_v^{-}:\mathfrak X_v^{-}\to X^{[d_v]}$$ is the forgetful morphism.

\subsection{Equal ranks and the proof of strange duality}
We now consider the case $r=s\geq 3$, when we simply have 
$$\Theta_w=c_v^{\star} L^{[d_v]}, \, \, \, \Theta_v =c_w^{\star} L^{[d_w]}.$$
Furthermore, tensor product gives a rational map defined away from codimension $2$ $$\tau^{+}:\mathfrak X_v^{+}\times \mathfrak X_w^{+}\dasharrow \mathfrak X^{+}\, \, \, (I_Z, z_B, I_T, t_B) \mapsto (I_Z \otimes I_{\widetilde T}, \, z_B \otimes t_B).$$ Here $$\mathfrak X^{+}=\{(U, u_B): a_B(U)=ru_B\}\subset X^{[d_v+d_w]}\times B.$$ 
In other words, $\mathfrak X^{+}$ is the fiber product
\begin{center}
$\xymatrix{ {{\mathfrak X}^{+}} \ar[d]^{\pi_2}  \ar[r]^{c} & X^{[d_v + d_w]} \ar[d]^{a_B}\\   B \ar[r]^{r} & B},$ 
\end{center}
where $c$ and $\pi_2$ are the natural projection maps, $$c ( U, u_B) = U, \, \, \, \pi_2 (U, u_B) = u_B.$$
The divisor $$\Theta^{+} \hookrightarrow \mathfrak X^{+}_v \times X^{+}_w$$ is the pullback under $\tau^+$ of the divisor  $$\theta^+=\{(U, u_B): h^0(I_U\otimes u_B \otimes L)\neq 0\},$$ with corresponding line 
bundle $${\mathcal O} \left ( \theta^{+} \right ) \simeq c^{\star} \, L^{[d_v + d_w]} \, \, \, \text{on} \, \, \, \mathfrak X^{+}.$$   

To identify the space of sections $H^0 \left (\mathfrak X^{+}, \, {\mathcal O} \left ( \theta^{+} \right ) \right)$, we let $B[r]$ denote the group of $r$-torsion points on $B$, and fix an isomorphism
$$r_{\star} {\mathcal O} \simeq \oplus_{\tau \in B[r]}\, \tau.$$
Then $$H^0\left (\mathfrak X^{+},  \, {\mathcal O} \left ( \theta^{+} \right ) \right)= H^0 (\mathfrak X^{+}, \, c^{\star} L^{[d_v +d_w]} )=H^0(X^{[d_v +d_w]}, \, L^{[d_v+d_w]}\otimes (c)_{\star} \mathcal O)=$$ $$ = \bigoplus_{\tau \in B[r]}  H^0(X^{[d_v+d_w]}, \, L^{[d_v+d_w]}\otimes a_B^{\star}\tau)=\bigoplus_{\tau \in B[r]} H^0(X^{[d_v+d_w]},\,  (L\otimes \tau)^{[d_v+d_w]}).$$

Under the isomorphism above, the divisor $\theta^{+}$ corresponds up to a ${\mathbb C}^{\star}$-scaling ambiguity to a tuple of sections, 
$$\theta^{+} \longleftrightarrow (s_\tau)_{\tau \in B[r]}, \, \, \, \, \, \, \, \, s_\tau \in H^0(X^{[d_v+d_w]},\,  (L\otimes \tau)^{[d_v+d_w]}).$$ The space of sections of $(L\otimes \tau)^{[d_v+d_w]}\to X^{[d_v+d_w]}$ can be identified with $$\Lambda^{d_v+d_w}H^0(L\otimes \tau),$$ which is one dimensional since $h^0(L\otimes \tau)=\chi(L\otimes \tau)=d_v+d_w.$ Furthermore, any non-zero section vanishes along the divisor 
$$\Theta_{L\otimes \tau} = \{I_U: \, \, \, h^0 (I_U \otimes L \otimes \tau) \neq 0 \}.$$ We have the important 
\begin{proposition}
\label{thetasplit}
For each $\tau, \, \, \, s_\tau$ is not the trivial section, hence it vanishes along $\Theta_{L\otimes \tau}.$
\end{proposition}

\noindent The proposition completes the proof of Theorem \ref{t1}. Indeed, as explained in \cite {tour}, each of the sections $s_{\tau}$ induces an isomorphism between spaces of sections 
$$\bigoplus_{\tau\in B[r]} s_{\tau}: \bigoplus_{\tau\in B[r]} H^0(X^{[d_v},\,  (L\otimes \tau)^{[d_v]})^{\vee}\to \bigoplus_{\tau\in B[r]} H^0(X^{[d_w]},\,  (L\otimes \tau)^{[d_w]}).$$ Since we argued that the strange duality map coincides with $\bigoplus_{\tau\in B[r]} s_{\tau}$, Theorem \ref{t1} follows. \footnote {Here, we assumed the determinant is $\mathcal O(\sigma+mf)$ where $f$ denotes the fiber over zero. If the determinant involves the fiber over a different point the same argument applies with the obvious changes in the definition of $L$.}

The proof of Theorem \ref{t2} is identical. \qed\vskip.1in

\noindent {\it Proof of the Proposition.} This is easily seen by restriction to  $X^{[d_v +d_w]}_F$, the fiber over zero of the addition map $$a_B: X^{[d_v +d_w]} \to B.$$ Letting $\mathfrak X^{+}_F$ be the fiber product
\begin{center}
$\xymatrix{ {{\mathfrak X}_F^{+}} \ar[d]^{\iota}  \ar[r]^{c} &  X^{[d_v + d_w]}_F \ar[d]^{\iota}\\  \mathfrak X^{+}  \ar[r]^{c} & X^{[d_v + d_w]}},$ 
\end{center}
we have a commutative diagram
\begin{center}
$\xymatrix{ { H^0 ({\mathfrak X}^{+}, \, \theta^{+}}) \ar[d]  \ar[r]  & \bigoplus_{\tau\in B[r]} H^0(X^{[d_v+d_w]},\,  (L\otimes \tau)^{[d_v+d_w]}) \ar[d] \\  H^0 (\mathfrak X^{+}_F, \theta^{+} )  \ar[r] & \bigoplus_{\tau\in B[r]} H^0(X^{[d_v+d_w]}_F,\,  (L\otimes \tau)^{[d_v+d_w]}) },$ 
\end{center}
where the horizontal maps are isomorphisms and the vertical maps are restrictions of sections.
The bottom isomorphism is particularly easy to understand since $${{\mathfrak X}_F^{+}}  \simeq  X^{[d_v + d_w]}_F \times B[r].$$
Let $\theta^{+}_F$ be the restriction of $\theta^{+}$ to $\mathfrak X^{+}_F$. Further restricting $\theta_F^+$ to $X_F^{[d_v+d_w]}\times \{\tau\}$ we obtain the divisor $\Theta_{L\otimes \tau}$. 
We claim that each of the divisors $\Theta_{L\otimes \tau}$ restricts nontrivially to $X_F^{[d_v+d_w]}$. This implies in turn that $s_{\tau}$ is not the trivial section for any $\tau$. The following lemma proves the claim above when $\tau=\mathcal O$; the arguments are identical for $\tau\neq \mathcal O$. \qed

\begin{lemma}\label{l5}
Consider the line bundle $L_{k,n} = \mathcal O(k\sigma+nf)$ on $X = B \times F,$ and assume that $k \geq 2$,  $n \geq 1.$ Then for a generic $I_Z \in X^{[kn]}_F,$ we have 
$$H^0 (I_Z \otimes L_{k,n}) =0.$$
\end{lemma}

\proof We argue by induction  on $n$.  For $n=1$, as established before, all sections of $L$ vanish along divisors of the form 
$$\bigsqcup_{i = 1}^{k} (B \times y_i) \cup f, \, \, \, \text{for} \, \, \, y_1 + \cdots + y_k = o_F.$$ A zero dimensional subscheme $I_Z$ in $X_F^{[k]}$ consisting of distinct points $z_i = (x_i, y_i)\in B\times F, \, \, 1\leq i \leq k$, must satisfy $$x_1 + \cdots + x_k = o_B.$$ 
Choose $Z$ so that $y_1 + \cdots + y_k \neq o_F,$ with $y_i\neq y_j$ for $i\neq j$, and so that $x_i \neq o_B$ for all $1\leq i \leq k.$ Then no section of $L$ vanishes at $Z$. 

To carry out the induction, note the exact sequence on $X$
$$0 \rightarrow  L_{k,n} \rightarrow  L_{k, n+1} \rightarrow {\mathcal O}_f (ko_F) \rightarrow 0,$$
and the associated exact sequence on global sections,
$$0 \rightarrow  H^0 (X, \, L_{k,n}) \rightarrow H^0 (X, \, L_{k, n+1}) \rightarrow H^0 (X, \, {\mathcal O}_f (ko_F)) \rightarrow 0.$$
Let $I_Z \in X_F^{[kn]}$ be so that $$H^0  (I_Z \otimes  L_{k,n}) = 0\, \, \, \text{and} \, \, \, Z \cap  \, (o_B \times F)  = \emptyset.$$ Then $H^1(I_Z\otimes L_{k, n})=0$. Choose $k$ additional points $$z_i = (o_B, y_i), \, 1\leq i \leq k, \, \, \, \text{so that} \, \, \, y_1 + \cdots y_k \neq o_F.$$
From the exact sequence on global sections, it follows that no section of $L_{k, n+1}$ vanishes at  $$Z \cup \{ z_1, \ldots, z_k\}\in X_F^{[k(n+1)]}.$$ Indeed, the exact sequence and the induction hypothesis on $Z$ show that $$H^0(X, \, L_{k, n+1}\otimes I_Z)\cong H^0(F,\,  \mathcal O_F(k o_F)\otimes I_Z) \cong H^0(F, \,\mathcal O_F(k o_F))$$ via restriction. But no section of ${\mathcal O}_F (k o_F)$ on the central fiber vanishes at the points $z_1, \ldots, z_k$, so no section of $L_{k, n+1}$ vanishes at $Z \cup \{ z_1, \ldots, z_k\}.$ We conclude that for a general $I_{{Z'}} \in X_F^{[k(n+1)]},$ we have $$H^0 (I_{{Z'}} \otimes L_{k,n+1}) =0.$$ This ends the proof of the lemma.\qed

\vskip.1in

\subsection {Proof of Theorem 2A} The strategy of proof is similar to that of Theorem \ref{t1}. We indicate the necessary changes. We write $$\det v=\mathcal O(\sigma+m_v f)\otimes Q_v^{-r}, \,\,\det w=\mathcal O(\sigma+m_w f)\otimes Q_w^{-s},$$ for line bundles $Q_v, Q_w$ of degree $0$ over $C$. We recall the birational isomorphism in Remark \ref{pa}: $$\mv\dasharrow \mathfrak X_v^+=\{(Z, z_C): z_C^r=\mathcal M_Z\}\subset X^{[d_v]}\times \text{Pic}^0(C)$$ where we set $$\mathcal M_Z=\mathcal O_C(a_C(Z)-d_v\cdot o).$$ This was noted when $Q_v$ is trivial in Section \ref{g2}, but the twist by $Q_v$ is an isomorphism of moduli spaces, yielding only a modified formula $$\fmd (V)=I_Z^{\vee}(-r\sigma+(\chi+\gbar)f)\otimes z_C^{-1}\otimes Q_v^{-1}.$$ There is a similar birational isomorphism $\mw\dasharrow \mathfrak X_w^+$.

The divisor $\Theta_{vw}^+\subset\mv\times \mw$ is identified with $$\Theta^+\subset \mathfrak X_v^+\times \mathfrak X_w^+$$ where $$\Theta^+=\{(Z, T, z_C, t_C): h^0(I_Z\otimes I_T\otimes z_C\otimes t_C\otimes L)\neq 0\},$$ for the line bundle $$L=\mathcal O((r+s)\sigma - (\chi_v+\chi_w+2\gbar) f)\otimes Q \implies \chi(L)=d_v+d_w.$$ Here, we wrote $$Q=Q_v\otimes Q_w$$ which by assumption is a generic line bundle of degree $0$ over the curve $C$. 
It can be seen that $L$ has no higher cohomology if $\chi_v+\chi_w\leq-3\gbar.$ In turn, when $r=s$, this is equivalent to the requirement $$d_v+d_w=-2r(\chi_v+\chi_w+\gbar)\geq 4r\gbar,$$ of the theorem. 

There are a few steps in the proof of Theorem 2A that need modifications from the genus $1$ case. They are: \begin {itemize} \item [(i)] the identification of the theta bundles. We carry this out for $r=s$ only, making use of the natural addition map $$\tau^+:\mathfrak X_v^+\times \mathfrak X_w^+\dasharrow \mathfrak X^+,$$ where $$ \mathfrak X^+=\{(U, u_C): \mathcal M_U=u_C^{r}\}\subset X^{[d_v+d_w]}\times \text{Pic}^0(C).$$ We have $\Theta^+=(\tau^+)^{\star} \theta^+$, for $$\theta^+=\{(U, u_C): h^0(I_U\otimes u_C\otimes L)\neq 0\}.$$ As before, we note the natural projections $$c:\mathfrak X^+\to X^{[d_v+d_w]},\,\,\, \text{pr}:\mathfrak X^+\to \text{Pic}^0(C).$$ We claim that \begin {equation}\label{pb}\mathcal O(\theta^+)=c^{\star} L^{[d_v+d_w]}\otimes \text{pr}^{\star} \mathcal P_{\alpha}^r,\end{equation} where $\mathcal P_{\alpha}$ is the line bundle over $\text{Pic}^0(C)$ associated to the point \begin{equation}\label{alp}\alpha=K_C(-2\gbar \cdot o)\otimes Q^{-2}\in \text{Pic}^0(C).\end{equation} Pulling back under $\tau^+$, it follows that $$\Theta_w=c_v^{\star} L^{[d_v]}\otimes \text{pr}^{\star}\mathcal P_{\alpha}^r,\,\, \Theta_v=c_w^{\star}L^{[d_w]} \otimes \text{pr}^{\star}\mathcal P_{\alpha}^r.$$ 

Before proving \eqref{pb}, we simplify notations. We write $$\ell =\chi(L)=d_v+d_w.$$ Over the Jacobian $$A=\text{Pic}^0(C),$$ we fix a principal polarization, for instance $$\Theta=\{y\in A: h^0(y\otimes \mathcal O_C(\bar g\cdot o))\neq 0\}.$$ The standardly normalized Poincar\'e bundle $$\mathcal P\to A\times A$$ takes the form $$\mathcal P=m^{\star} \Theta^{-1}\otimes (\Theta \boxtimes \Theta).$$ This differs up to a sign from the usual conventions, but it is compatible with our conventions on the Abel-Jacobi embedding $$C\to A,\,\, x\to \mathcal O_C(x-o),$$ in the sense that that $\mathcal P|_{A\times \alpha}$ restricts to $\alpha$ over $C$. For degree $\ell$, we consider  the Abel-Jacobi map $$\pi:C^{(\ell)}\to A, \,\, D\mapsto \mathcal O(D-\ell \cdot o),$$ and set  $$\mathcal P_C^{(\ell)}=(\pi,1)^{\star} \mathcal P\to C^{(\ell)}\times A.$$ The bundle $\mathcal P^{(\ell)}$ satisfies $$\mathcal P^{(\ell)}|_{C^{(\ell)}\times \{y\}}\cong y^{(\ell)},\,\, \mathcal P^{(\ell)}|_{\ell [o]\times A} \text{ is trivial over} A.$$ When $\ell=1$, $\mathcal P^{(1)}$ is the Poincar\'e bundle $$\mathcal P_C\to C\times A$$ normalized over $o$.

With these notations out of the way, we first consider $\mathcal O(\theta^+)$ over the product  $X^{[\ell]}\times A.$ We claim that \begin {equation}\label{h1} \mathcal O(\theta^+)=c^{\star} L^{[\ell]}\otimes (a_C, 1)^{\star} \mathcal P^{(\ell)}\otimes \text{pr}^{\star}_A \mathcal M\end{equation} for some line bundle $\mathcal M\to A.$ This follows from the see-saw theorem. The restriction of $\mathcal O(\theta^+)$ to $X^{[\ell]}\times \{u_C\}$ is $$(L\otimes u_C)^{[\ell]}=L^{[\ell]}\otimes a_C^{\star} u_C^{(\ell)}.$$ This agrees with the restriction of $$c^{\star} L^{[\ell]}\otimes (a_C, 1)^{\star} \mathcal P^{(\ell)}$$ and establishes \eqref{h1}. To identify $\mathcal M$, we restrict to $\{U\}\times A,$ where $U$ is a length $\ell$ subscheme of $X$ supported over $o\times o_F$. We obtain $$\mathcal M=\det \rp p_{!} (\mathcal P_C\otimes q^{\star}(L\otimes I_U))^{-1}.$$ We can rewrite $\mathcal M$ expressing in $K$-theory $$I_U=\mathcal O-\ell\cdot \mathcal O_{o\times o_F}.$$ Recalling the normalization of $\mathcal P_C$ over $o$, and that $$L=\left(\mathcal O_C(t\cdot o)\otimes Q\right)\boxtimes \mathcal O_F(2r\cdot o_F),$$ for $t=-(\chi_v+\chi_w+2\gbar)$, we obtain 
$$\mathcal M=\det \rp p_{!}(\mathcal P_C\otimes q^{\star} L)^{-1}=\det \rp p_{!} (\mathcal P_C\otimes q^{\star} (\mathcal O_C(t\cdot o)\otimes Q))^{-2r}\\$$$$=\det \rp p_! (\mathcal P_C\otimes q^{\star} (\mathcal O_C(\gbar \cdot o)\otimes Q))^{-2r}=\Theta^{2r}\otimes \mathcal P_{-Q}^{2r}.$$
Therefore $$\mathcal O(\theta^+)=c^{\star} L^{[\ell]}\otimes (a_C, 1)^{\star}\circ (\pi, 1)^{\star}\mathcal P \otimes \text{pr}_A^{\star} (\Theta\otimes \mathcal P_{-Q})^{2r}.$$  Over $\mathfrak X^+\hookrightarrow X^{[\ell]}\times A $, we have the following commutative diagram 
\begin{center}
$\xymatrix{ {{\mathfrak X}^{+}}\ar[d]^{\text{pr}_A}  \ar[r]^{(a_C,1)} & C^{(\ell)}\times A \ar[d]^{(\pi, 1)}\\   A \ar[r]^{(r,1)} & A\times A}$ \end{center}
so that $$(\pi, 1)\circ (a_C, 1)= (r, 1)\circ \text{pr}_A.$$ Hence, we obtain \begin{eqnarray*}\mathcal O(\theta^+)&=&c^{\star} L^{[\ell]}\otimes\text{pr}_{A}^{\star}((r, 1)^{\star} \mathcal P\otimes \Theta^{2r}\otimes \mathcal P_{-Q}^{2r})\\ &=&c^{\star}L^{[\ell]}\otimes \text{pr}_{A}^{\star} (r^{\star}\Theta\otimes\Theta\otimes (r+1)^{\star}\Theta^{-1}\otimes \Theta^{2r}\otimes \mathcal P_{-Q}^{2r})\\&=&c^{\star} L^{[\ell]}\otimes \text{pr}_A^{\star}(\Theta\otimes (-1)^{\star}\Theta^{-1}\otimes \mathcal P_{-Q}^2)^{r}=c^{\star}L^{[\ell]}\otimes \text{pr}_A^{\star}\mathcal P^r_\alpha,\end{eqnarray*} as claimed.

\item [(ii)] Next, we identify the space of sections $H^0(\mathfrak X^+, \mathcal O(\theta^+))$ using the cartesian diagram \begin{center}
$\xymatrix{ {{\mathfrak X}^{+}} \ar[d]^{\text{pr}_A}  \ar[r]^{c} & X^{[\ell]} \ar[d]^{f}\\   A \ar[r]^{r} & A},$ 
\end{center} where $f=\pi\circ a_C$. 
We have \begin{eqnarray*}H^0(\mathfrak X^+, \mathcal O(\theta^+))&=&
H^0(\mathfrak X^+, c^{\star} L^{[\ell]}\otimes \text{pr}_A^{\star} \mathcal P_{\alpha}^r)=H^0(X^{[\ell]}, L^{[\ell]}\otimes c_{\star}\text{pr}_A^{\star} \mathcal P_{\alpha}^r)\\&=&H^0(X^{[\ell]}, L^{[\ell]}\otimes f^{\star} r_{\star} \mathcal P_{\alpha}^r)=H^0(X^{[\ell]}, L^{[\ell]}\otimes f^{\star} r_{\star} r^{\star} \mathcal P_{\alpha})\\&=&\bigoplus_{\tau\in A[r]} H^0(X^{[\ell]}, L^{[\ell]}\otimes f^{\star}(\mathcal P_{\alpha}\otimes \tau))\\&=&\bigoplus_{\tau\in A[r]} H^0(X^{[\ell]}, L^{[\ell]}\otimes a_C^{\star}(\alpha\otimes \tau))\\&=&\bigoplus_{\tau\in A[r]} H^0(X^{[\ell]}, (L\otimes \alpha\otimes \tau)^{[\ell]}).\end{eqnarray*} Similar expressions hold over each of the factors $\mathfrak X_v^+$ and $\mathfrak X_w^+$. 
\item [(iii)] Finally, we need a suitable analogue of Lemma \ref{l5}. This concerns subschemes $Z$ in $X$, belonging to $$X_F^{[\ell]}=\{Z: a_C(Z) \text{ is rationally equivalent to } \ell \cdot o\}.$$ We show that if $$L_{k,n}=\mathcal O(k\sigma+nf)$$ for $k\geq 2g, n\geq g,$ then for all $r$-torsion line bundles $\tau$ over $C$, we have $$H^0(L_{k, n}\otimes \alpha \otimes \tau \otimes I_Z)=0$$ provided $Z$ is generic in $X_F^{[k(n-\gbar)]}$, and $Q$, which appears in the definition of $\alpha$ in \eqref{alp}, is generic over $C$. This is then applied to the line bundle $L$ appearing in the expression of the theta bundles above. 

We induct on $n$, starting with the base case $n=g$. Just as in genus $1$, for generic $Q$, the sections of $L_{k, g}\otimes \alpha \otimes \tau$ vanish along divisors of the form $$\bigsqcup_{i = 1}^{k} (C \times y_i) \cup f_{p_1}+\ldots+f_{p_g}, \, \, \, \text{for} \, \, \, y_1 + \cdots + y_k = o_F,$$ and some $p_1, \ldots, p_g\in C$. 
Indeed, it suffices to explain that $$h^0(\mathcal O_C(g\cdot o)\otimes \alpha \otimes \tau)=1,$$ the unique section vanishing at $g$ points $p_1(o), \ldots, p_g(o),$ which depend on $o$. Equivalently, recalling the definition of $\alpha$ in \eqref{alp}, and using Riemann-Roch and Serre-duality, 
we prove that $$h^0(\tau^{-1}\otimes \mathcal O_C((g-2)\cdot o)\otimes Q^2)=0.$$ As the Euler characteristic of the above line bundle is $-1$, the space of sections vanishes for generic $Q$ of degree $0$. 

To complete the proof of the base case, fix a generic $Q$ as above. It suffices to explain that for all $k\geq 2g$, we can find $x_1, \ldots, x_{k}\in C$ such that $$x_1+\ldots+x_k\equiv k\cdot o, \, x_i\neq p_j \text{ for all } i, j.$$ Choosing a section from the non-empty set $$H^0(\mathcal O_C(k\cdot o))-\bigcup_{j=1}^{g} H^0(\mathcal O_C(k\cdot o-p_j))-H^0(\mathcal O_C((k-1)\cdot o)),$$ we let $x_1, \ldots, x_k$ be its zeros. This ends the argument for the base case. 

The inductive step does not require any changes from the original Lemma \ref{l5}.  
\end {itemize}

The proof of Theorem 2A is now completed. \qed

\end {document}